\documentclass[onecolumn,final,a4paper,10pt,sort&compress]{elsarticle}

\usepackage[a4paper]{geometry}
\usepackage{amsfonts}
\usepackage{amsmath}
\usepackage[english]{babel}
\usepackage{bbm}
\usepackage{color}
\usepackage{graphicx}
\newcommand{\url}[1]{\texttt{#1}}
\usepackage[utf8x]{inputenc}
\usepackage{lmodern}
\usepackage{multicol}

\newcommand{\smalleq}[1]{
    \begin{small}\begin{align*}
    #1
    \end{align*}\end{small}}
    
\newcommand{\GAs}[1] {\textbf{G{#1}$^{\rm{2D}}$}}
\newcommand{\GAv}[1] {\textbf{G{#1}}}

\newcommand{\rhos}[1] {\varrho^{\rm{2D}}_{#1}}
\newcommand{\rhov}[1] {\varrho_{#1}}

\newcommand{\Dtetuni}{\calD_{\rm{tet-uniform}}}
\newcommand{\Dtetani}{\calD_{\rm{tet-anisotropic}}}
\newcommand{\Dtetpar}{\calD_{\rm{tet-parallel}}}
\newcommand{\Dtetbcl}{\calD_{\rm{tet-bcl}}}
\newcommand{\Dtetpoisson}{\calD_{\rm{tet-poisson}}}
\newcommand{\Dtetrandom}{\calD_{\rm{tet-random}}}

\newcommand{\Dhexuni}{\calD_{\rm{hex-uniform}}}
\newcommand{\Dhexani}{\calD_{\rm{hex-anisotropic}}}
\newcommand{\Dhexpar}{\calD_{\rm{hex-parallel}}}

\newcommand{\Dvorobcl}{\calD_{\rm{voro-bcl}}}
\newcommand{\Dvoropoisson}{\calD_{\rm{voro-poisson}}}
\newcommand{\Dvororandom}{\calD_{\rm{voro-random}}}

\newcommand{\Dpolypar}{\calD_{\rm{poly-parallel}}}
\newcommand{\Dpolypoisson}{\calD_{\rm{poly-poisson}}}
\newcommand{\Dpolyrandom}{\calD_{\rm{poly-random}}}

\newcommand{\rP}{r_\P}
\newcommand{\rF}{r_\F}

\newcommand{\sqbra}[1]{\left[{#1}\right]}
\newcommand{\cubra}[1]{\left\{{#1}\right\}}
\newcommand{\REAL}{\mathbbm{R}}

\newcommand{\calD}{\mathcal{D}}
\newcommand{\calE}{\mathcal{E}}
\newcommand{\calF}{\mathcal{F}}

\newcommand{\calH}{\mathcal{H}}
\newcommand{\calI}{\mathcal{I}}

\newcommand{\calT}{\mathcal{T}}

\newcommand{\calV}{\mathcal{V}}

\newcommand{\as}{a}

\newcommand{\fs}{f}
\newcommand{\gs}{g}

\newcommand{\ms}{m}

\newcommand{\qs}{q}

\newcommand{\us}{u}
\newcommand{\vs}{v}

\newcommand{\xs}{x}
\newcommand{\ys}{y}
\newcommand{\zs}{z}

\newcommand{\Bs}{B}
\newcommand{\Cs}{C}

\newcommand{\Ns}{N}

\newcommand{\Vs}{V}

\newcommand{\xsP}{\xs_{\P}}
\newcommand{\ysP}{\ys_{\P}}
\newcommand{\zsP}{\zs_{\P}}

\newcommand{\PS}[1]{\mathbbm{P}_{#1}}

\newcommand{\HONE}   {H^1}

\newcommand{\LTWO}  {L^2}
\newcommand{\LINF}  {L^{\infty}}

\newcommand{\LS}[1] {L^{#1}}
\newcommand{\HS}[1] {H^{#1}}

\newcommand{\CS}[1] {C^{#1}}

\renewcommand{\P} {E}
\newcommand  {\F} {\textsf{f}}            
\newcommand  {\E} {\textsf{e}}            
\newcommand  {\V} {\textsf{v}}            

\newcommand{\hh}{h}
\newcommand{\Th}{\Omega_{\hh}}

\newcommand{\xv}{\mathbf{x}}
\newcommand{\yv}{\mathbf{y}}

\newcommand{\xvP}{\xv_{\P}}        
\newcommand{\xvF}{\xv_{\F}}        
\newcommand{\xvV}{\xv_{\V}}        

\newcommand{\hP}{\hh_{\P}}
\newcommand{\hF}{\hh_{\F}}
\newcommand{\hE}{\hh_{\E}}

\newcommand{\mP}{\ABS{\P}}
\newcommand{\mF}{\ABS{\F}}

\newcommand{\dV}{\,dV}

\newcommand{\fsh}{\fs_{\hh}}

\newcommand{\ush}{\us_{\hh}}

\newcommand{\vsh}{\vs_{\hh}}

\newcommand{\Bsh}{\Bs^{\hh}}

\newcommand{\asP}{\as^{\P}}

\newcommand{\ash}{\as_{\hh}}

\newcommand{\ashP}{\as^{\P}_{\hh}}

\newcommand{\SP} {S^{\P}}

\newcommand{\nlen}{\hspace{-0.2mm}}

\newcommand{\snorm}  [2]{\vert#1\vert_{#2}}

\newcommand{\norm}   [2]{|\nlen|#1|\nlen|_{#2}}

\newcommand{\ABS}    [1]{\vert#1\vert}

\newcommand{\Vsh}{V^{\hh}_{}}
\newcommand{\Vsht}{\widetilde{V}^{\hh}_{}}

\newcommand{\Piz}[1]{\Pi^{0}_{#1}}
\newcommand{\PinP}[1]{\Pi^{\nabla,\P}_{#1}}
\newcommand{\PinF}[1]{\Pi^{\nabla,\F}_{#1}}
\newcommand{\PizP}[1]{\Pi^{0,\P}_{#1}}

\newcommand{\restrict}[2]{{#1}_{\vert{#2}}}

\newcommand{\EOD}{\end{document}}

\newcommand{\bil}[2]{\langle#1,#2\rangle}

\newcommand{\bP} {\partial \P}  

\newcommand{\Vsg}{\Vs_{\gs}}
\newcommand{\Vszr}{\Vs_{0}}

\newcommand{\Vshg}{V^{\hh}_{\gs}}
\newcommand{\Vshzr}{V^{\hh}_{0}}

\newcommand{\gsI}{\gs_{I}}

\begin{document}

\begin{frontmatter}


  \title{Polyhedral Mesh Quality Indicator for the Virtual Element Method}
  
  
  \author[IMATI]{T. Sorgente\corref{cor1}}
  \ead{tommaso.sorgente@ge.imati.cnr.it}
  
  \author[IMATI]{S. Biasotti}
  \ead{silvia.biasotti@ge.imati.cnr.it}

  \author[IMATI] {G. Manzini}
  \ead{marco.manzini@imati.cnr.it}
  
  \author[IMATI]{M. Spagnuolo}
  \ead{michela.spagnuolo@ge.imati.cnr.it}

  \cortext[cor1]{Corresponding authors}


  \address[IMATI]{ Istituto di Matematica Applicata e Tecnologie
    Informatiche, Consiglio Nazionale delle Ricerche, Italy}


  \begin{abstract}
    We present the design of a \emph{mesh quality indicator} that can
    predict the behavior of the Virtual Element Method (VEM) on a given mesh
    family or finite sequence of polyhedral meshes (dataset).
    The mesh quality indicator is designed to measure the violation of
    the mesh regularity assumptions that are normally considered in
    the convergence analysis.
    We investigate the behavior of this new mathematical tool on
    the lowest-order conforming approximation of the
    three-dimensional Poisson equation.
    This work also assesses the convergence rate of the VEM when
    applied to very general polyhedral meshes, including non
    convex and skewed three-dimensional elements. 
    Such meshes are created within an original mesh generation framework, which is designed to allow the generation of meshes with very different sizes, connectivity and geometrical properties.
    The obtained results show a significant correlation between the quality measured a priori by the indicator and the effective performance of the VEM.
  \end{abstract}
  
  
  \begin{keyword}
    virtual element method,
    polyhedral mesh,
    mesh regularity assumptions,
    mesh quality indicators,
    small edges,
    small faces,
    3D Poisson problem,
    optimal convergence
    \\~\\
    \textbf{AMS subject classification:} 65N12; 65N15
  \end{keyword}
  
\end{frontmatter}
  

\raggedbottom
\setcounter{secnumdepth}{4}
\setcounter{tocdepth}{4}


\section{Introduction}
\label{sec:intro}
The Virtual Element Method
(VEM)~\cite{BeiraodaVeiga-Brezzi-Cangiani-Manzini-Marini-Russo:2013}
is a Galerkin projection method such as the Finite Element Method
(FEM).
The major difference between VEM and FEM is that VEM does not require
the explicit evaluation of the basis functions and their gradients,
which are integrated in the variational formulation.
In the VEM, the basis functions are formally defined as the solutions
to suitable partial differential equation problems formulated in every
mesh element, and they are dubbed as \emph{virtual} since they are never
explicitly evaluated.
The method relies on some special polynomial projections of the basis
functions and their derivatives, which are computable from a careful
choice of the degrees of freedom.
Using these projection operators, we define the bilinear forms and
linear functionals of the variational formulation.
The resulting schemes can be proved to be consistent with polynomials
of a given degree and this property determines the accuracy of the
discretization, while the stability of the method, which implies its
well-posedness, is ensured by introducing in the formulation a
suitable stabilization term.

This computational approach is extremely powerful and offers indeed
several potential advantages with respect to the FEM.
In fact, we can easily build approximation spaces that work on very
general meshes, including meshes whose elements are generic-shaped polygons in 2D
and polyhedra in 3D (also called \emph{polytopes} for short when we
do not need to specify the number of dimensions).
These finite dimensional spaces can show arbitrary global regularity,
a property useful in the discretization of high order partial
differential equations such as, for example, the polyharmonic
problems~\cite{Antonietti-Manzini-Verani:2019:CAMWA:journal}.
The VEM has thus been proved to be very successful, and an incomplete
list of significant applications on general meshes includes, for
example, the works of Refs.~\cite{%
  Antonietti-Manzini-Verani:2018,%
  Brezzi-Marini:2013,%
  BeiraodaVeiga-Manzini:2014,%
  BeiraodaVeiga-Manzini:2015,%
  BeiraodaVeiga-Mora-Vacca:2019,%
  BeiraodaVeiga-Manzini-Mascotto:2019,%
  BeiraodaVeiga-Dassi-Manzini-Mascotto:2021,%
  Benedetto-Berrone-Pieraccini-Scialo:2014,%
  Berrone-Borio-Manzini:2018:CMAME:journal,%
  Benvenuti-Chiozzi-Manzini-Sukumar:2019:CMAME:journal,%
  Cangiani-Gyrya-Manzini:2016,%
  Cangiani-Georgoulis-Pryer-Sutton:2016,%
  Cangiani-Manzini-Sutton:2017,%
  Certik-Gardini-Manzini-Vacca:2018:ApplMath:journal,
  Certik-Gardini-Manzini-Mascotto-Vacca:2019:CAMWA:journal,%
  Gardini-Manzini-Vacca:2019:M2AN:journal,%
  Mora-Rivera-Rodriguez:2015,%
  Natarajan-Bordas-Ooi:2015,%
  Paulino-Gain:2015,%
  Perugia-Pietra-Russo:2016,%
  Wriggers-Rust-Reddy:2016}.
A detailed description of the state of the art can also be found in
the very recent collection of thematic
articles~\cite{SEMA-SIMAI-BOOK:2022}.

A major property of the VEM is that this method is especially suited to
solving partial differential equations (PDEs) on polygonal and
polyhedral meshes.
The VEM is indeed highly versatile in the admissible meshes.
This property is very useful when we consider a mesh adaptation
strategy that allows the mesh to be locally refined in those parts of
the domain requiring greater accuracy in order to improve the
numerical approximation.
A second major property of the method, which was noted form the very
beginning of the VEM history, is its extreme robustness with respect
to mesh deformations.
For example, the VEM can be used for problems with oddly shaped
material interfaces to which the mesh must be conformal, or where the
boundary of the computational domain deforms in time.

The convergence of the VEM can be proved under different sets of mesh
regularity assumptions, which impose restrictions on the meshes that
we use in practical calculations.
A detailed review of these results for the numerical approximation of
the two-dimensional Poisson equation can be found in the book
chapter~\cite{sorgente2021vem}.
In our first work~\cite{sorgente2021role}, we investigated this aspect
in a very extensive way and we noted that even relaxing part of these
assumptions does not seem to impact significantly on the optimal
convergence behavior of the VEM when applied to the discretization of
the two-dimensional Poisson equation.
Such an approximation was shown to be robust and accurate even on
meshes with highly irregular structures or meshes with highly skewed
elements, as for example, those ones that could be outcome from
mesh adaptation algorithms and refinement.

We also designed a \emph{mesh quality indicator}, which is a
mathematical tool that predicts the behavior of the virtual element
method on a given sequence of refined meshes \emph{without solving the
differential problem}.
This indicator can be useful to select the kind of meshes to which we
can apply the VEM and, also, to improve the quality of the meshes that
are outcome by some adaptation or agglomeration algorithm.

The present work aims at generalizing such indicator to
\emph{three-dimensional polytopal meshes}.
Again, the indicator is useful to predict the behavior of the VEM over a given sequence of meshes \emph{before} applying the numerical method itself, so that we can evaluate if the given
mesh sequence is suited to the virtual element discretization.
The paper starts from an overview of the geometrical assumptions to guarantee the convergence
introduced in the literature of the conforming VEM for the Poisson
equation in three dimensions.
Similarly to what has been done in \cite{sorgente2021role}, from these
assumptions we derive our \textit{quality indicator} for polyhedral
meshes, which uniquely depends on the geometry of the mesh elements
and measures the violation of the geometrical assumptions.
Then, we define a mesh generation framework, consisting in a number of
\textit{sampling strategies} and \textit{meshing techniques}, which
allows us to build sequences of polyhedral meshes (datasets) with
different connectivity.
Many of the so-generated datasets violate the considered geometrical
assumptions, and this fact allows a correlation analysis between
such assumptions and the VEM performance.
We are then ready to measure the quality of each dataset and solve the
Poisson problem over it with the lowest-order conforming VEM, which is
built upon the elemental virtual element spaces that contain the
linear polynomials as a subspace.
We focus on the linear virtual element method for computational
simplicity and because we think that this is the most interesting
choice in three-dimensional practical engineering applications,
leaving a similar investigation for the higher-order formulations to
future work.
We experimentally show how the VEM presents a good convergence rate in
the majority of the cases, underperforming only in very few situations.
We also show a strict correspondence between the values measured by the quality indicator and
the performance of the VEM on a given mesh, or dataset, both in terms
of approximation error and convergence rate.

\subsection{Notation and technicalities}
\label{subsec:notation}
We use the standard definition and notation of Sobolev spaces, norms
and seminorms, cf.~\cite{Adams-Fournier:2003}.
Let $k$ be a nonnegative integer number.
The Sobolev space $\HS{k}(\omega)$ consists of all square integrable
functions with all square integrable weak derivatives up to order $k$
that are defined on the open, bounded, connected subset $\omega$ of
$\REAL^{d}$, $d=1,2$.
As usual, if $k=0$, we prefer the notation $\LTWO(\omega)$.
Norm and seminorm in $\HS{k}(\omega)$ are denoted by
$\norm{\cdot}{k,\omega}$ and $\snorm{\cdot}{k,\omega}$, while for the
inner product in $\LTWO(\omega)$ we prefer the integral notation.
We denote the space of polynomials of degree less than or equal to
$k\geq0$ on $\omega$ by $\PS{k}(\omega)$ and conventionally assume
that $\PS{-1}(\omega)=\{0\}$.
In our implementation, we consider the orthogonal basis on every mesh
edge through the univariate Legendre polynomials and inside every mesh
cell provided by the Gram-Schmidt algorithm applied to the standard
monomial basis.

\subsection{Outline}
The paper is organized as follows.
In Section~\ref{sec:vem}, we present the VEM in three dimensions and
the convergence results for the Poisson equation with Dirichlet
boundary conditions.
In Section~\ref{sec:quality}, we report the geometrical assumptions
on the mesh elements that are used in the literature to guarantee the
convergence of the VEM, both in 2D and in 3D.
We then introduce the mesh quality indicator for 3D polyhedral meshes, which extends the indicator presented in \cite{sorgente2021role} for polygonal meshes.
In Section~\ref{sec:datasets}, we define and build a number of polyhedral
datasets and analyze their geometrical properties.
In Section~\ref{sec:results}, we apply the mesh quality indicator to
the datasets and solve the problem with the VEM, looking for
correlations between the quality and the approximation errors.
In Section~\ref{sec:conclusions}, we offer our final remarks and
discuss future developments and work.
All the meshes used in this work are available at\\
\centerline{\url{https://github.com/TommasoSorgente/vem-indicator-3D-dataset}}
for download.

\section{The Virtual Element Method on polyhedral meshes}
\label{sec:vem}

In this section, we briefly review the lowest-order virtual element
method in three space dimensions for the Poisson equation in primal
form.
A more detailed presentation of these concepts can be found in
\cite{%
  BeiraodaVeiga-Brezzi-Cangiani-Manzini-Marini-Russo:2013,%
  BeiraodaVeiga-Brezzi-Marini-Russo:2014,%
  Ahmad-Alsaedi-Brezzi-Marini-Russo:2013}
where the VEM with arbitrary order of accuracy is presented.
The extension to a more general second-order three-dimensional elliptic
problem including a reaction term is found
in~\cite{da2017high}.

\subsection{Mesh generalities}
\label{subsec:mesh:generalities}
The virtual element method is formulated on the mesh family $\calT=$
$\big\{\Th\big\}_{h}$, where each mesh $\Th$ is a partition of the
computational domain $\Omega$ into nonoverlapping polyhedral elements
$\P$ and is labeled by the mesh size parameter $\hh$ that is defined
below.
A polyhedral element $\P$ is a compact subset of $\REAL^3$ with
boundary $\partial\P$, volume $\mP$, barycenter
  $\xvP=(\xsP,\ysP,\zsP)^T$, and diameter
$\hP=\sup_{\xv,\yv\in\P}\vert\xv-\yv\vert$.
The set of the mesh elements $\Th$ form a finite cover of $\Omega$
such that $\overline{\Omega}=\cup_{\P\in\Th}\P$ and the mesh size
labeling each mesh $\Th$ is defined by $\hh=\max_{\P\in\Th}\hP$.
We assume that the mesh sizes of the mesh family $\calT$ are in a
countable subset $\calH$ of the real line $(0, +\infty)$ having $0$ as
its unique accumulation point.
A mesh face $\F$ is a planar, two-dimensional subset of $\REAL^3$ with
area $\mF$, barycenter $\xvF$ and diameter
$\hF=\sup_{\xv,\yv\in\F}\vert\xv-\yv\vert$ and we denote the set of
mesh faces by $\calF_{\hh}$.
We shall consider on every mesh face a local coordinate system
$(\xi,\eta)$.
A mesh edge $\E$ is a straight one-dimensional subset of $\REAL^3$
with length $\hE$ and we denote the set of mesh edges by
$\calE_{\hh}$.
We shall consider on every mesh edge a local coordinate system $s$.
A mesh vertex $\V$ has three-dimensional position vector $\xvV$ and we
denote the set of mesh vertices by $\calV_{\hh}$.

\medskip
According to the notation introduced in Section~\ref{subsec:notation},
we denote the space of linear polyomials defined on $\P$, $\F$, and
$\E$ by $\PS{1}(\P)$, $\PS{1}(\F)$ and $\PS{1}(\E)$, respectively, and
the space of piecewise linear polynomials on the whole mesh $\Th$ by
$\PS{1}(\Th)$.
Accordingly, if $\qs\in\PS{1}(\Th)$ then it holds that
$\restrict{\qs}{\P}\in\PS{1}(\P)$ for all $\P\in\Th$.

In the implementation, we consider the following basis of $\PS{1}(\P)$
in each element $\P$:
\begin{align*}
  \ms^{\rm{3D}}_0(\xs,\ys,\zs) &= 1,\quad
  \ms^{\rm{3D}}_1(\xs,\ys,\zs) = \frac{\xs-\xsP}{\hP},\quad
  \ms^{\rm{3D}}_2(\xs,\ys,\zs) = \frac{\ys-\ysP}{\hP},\\[0.5em]
  \ms^{\rm{3D}}_3(\xs,\ys,\zs) &= \frac{\zs-\zsP}{\hP}.
\end{align*}
Similarly, we consider the following basis of $\PS{1}(\F)$ on each
face $\F$:
\begin{align*}
  \ms^{\rm{2D}}_0(\xi,\eta) = 1,\quad
  \ms^{\rm{2D}}_1(\xi,\eta) = \frac{\xi -\xi_{\F}}{\hF},\quad
  \ms^{\rm{2D}}_2(\xi,\eta) = \frac{\eta-\eta_{\F}}{\hF},
\end{align*}
where we recall that $(\xi,\eta)$ is the local coordinate system
defined on $\F$ and $\xvF=(\xi_{\F},\eta_{\F})^T$ is the center of the
face.
We also use the following basis of $\PS{1}(\E)$ on every edge $\E$:
\begin{align*}
  \ms^{\rm{1D}}_0(s) = 1,\quad
  \ms^{\rm{1D}}_1(s) = \frac{s -s_{\E}}{\hE},
\end{align*}
where we recall that $(\xi,\eta)$ is the local coordinate system
defined on $\E$ and $s_{\P}$ is the center (mid-point) of the edge.

\subsection{The model problem}
\label{subsec:model:problem}
Let $\Omega\subset\REAL^3$ be an open, bounded, simply connected,
convex domain with Lipschitz boundary $\Gamma$.
For exposition's sake, we assume that $\Omega$ is a polyhedral domain
and its boundary $\Gamma$ is given by the union of a subset of the
faces in $\calF$.
We consider the diffusion problem
\begin{align}
  -\Delta\us &= \fs\phantom{\gs} \quad\textrm{in}~\Omega,\label{eq:Poisson:a}\\[0.5em]
  \us        &= \gs\phantom{\gs} \quad\textrm{on}~\Gamma,\label{eq:Poisson:b}
\end{align}
where $\fs\in\LS{2}(\Omega)$ is the load term and
$\gs\in\HS{\frac12}(\Gamma)$ the Dirichlet boundary data.
The variational form of this problem reads as:
\begin{align}
  \mbox{\textit{Find $\us\in\Vsg$ such that}}\quad
  \int_{\Omega}\nabla\us\cdot\nabla\vs = \int_{\Omega}\fs\vs \quad\forall\vs\in\Vszr,
  \label{eq:varform}
\end{align}
where $\Vsg=\{ \vs\in\HS{1}(\Omega):\,\restrict{\vs}{\Gamma}=\gs \}$
and $\Vszr=\{ \vs\in\HS{1}(\Omega):\,\restrict{\vs}{\Gamma}=0 \}$.
The well-posedness of this problem follows from the coercivity and
continuity of the bilinear form on the left-hand side
of~\eqref{eq:varform} and the boundedness on the linear functional of
the right-hand side of~\eqref{eq:varform} and can be proved by an
application of the Lax-Milgram theorem,
see~\cite[Theorem~2.7.7]{Brenner-Scott:2008}.

\subsection{Virtual elements on polygonal faces}
\label{subsec:VEM:2D:faces}
The three-dimensional conforming virtual element space is built
recursively on top of the virtual element spaces defined on the
polyhedral faces.
Under the assumptions presented in
section~\ref{subsec:mesh:generalities}, each polyhedral face $\F$ is a
two-dimensional planar polygon.
In order to define the virtual element space on $\F$, we preliminarly
introduce the virtual element space
\begin{align}
  \Vsht(\F) := \Big\{
  \vsh\in\HS{1}(\F)\cap\CS{0}(\F):\,
  \restrict{\vsh}{\E}\in\PS{1}(\E)\,\forall\E\in\partial\F,\,
  \Delta\vsh\in\PS{1}(\F)
  \Big\}.
  \label{eq:VEM:space:2D:extended}
\end{align}
Next, we consider the linear, bounded operators
$\lambda_{\ell}:\Vsh(\F)\to\REAL$ such that
$\lambda_{\ell}(\vsh)=\vsh(\V_{\ell})$ is the value at the $\ell$-th
vertex of $\partial\F$ and the elliptic projection operator
$\PinF{\F}:\Vsh(\F)\to\PS{1}(\F)$, which is such that
\begin{align}
  \int_{\F}\nabla\big(\vsh-\PinF{1}\vsh\big)\cdot\nabla\qs &= 0\quad\forall\qs\in\PS{1}(\F),\label{eq:PinF:def:A}\\[0.5em]
  \int_{\partial\F}\big(\vsh-\PinF{1}\vsh\big)               &= 0.                             \label{eq:PinF:def:B}
\end{align}
An integration by parts shows that we can reduce the integral
in~\eqref{eq:PinF:def:A} to a boundary integral that can be evaluated
by a splitting on the elemental edges.
On each edge we apply the trapezoidal rule, which only depends on the
vertex values, i.e., the values $\lambda_{\ell}(\vsh)$, and returns
the exact integral value since the trace of $\vsh$ is linear,
cf.~\cite{BeiraodaVeiga-Brezzi-Cangiani-Manzini-Marini-Russo:2013,BeiraodaVeiga-Brezzi-Marini-Russo:2014}.
The integral in~\eqref{eq:PinF:def:B} is evaluated in the same way and
also depends only on the vertex values of $\vsh$.

Then, we define the two-dimensional virtual element space on $\F$ as
\begin{align}
  \Vsh(\F) := \bigg\{
  \vsh\in\Vsht(\F):\,\int_{\F}\vsh\qs=\int_{\F}\big(\PinP{1}\vsh\big)\qs
  \quad\forall\qs\in\PS{1}(\P)
  \bigg\}.
  \label{eq:VEM:space:2D:enhanced}
\end{align}
A major property of definition~\eqref{eq:VEM:space:2D:enhanced} is
that the linear polynomials over $\F$ are a subspace of $\Vsh(\F)$;
formally, we can write that
$\PS{1}(\F)\subseteq\Vsh(\F)\subseteq\Vsht(\F)$.
Moreover, the virtual element functions $\vsh\in\Vsh(\F)$ are uniquely
identified by their vertex values. i.e., $\lambda_{\ell}(\vsh)$.
A proof of the unisolvence of the set $\{\lambda_{\ell}(\vsh)\}$ in
$\Vsh(\F)$ can be
found in~\cite{Ahmad-Alsaedi-Brezzi-Marini-Russo:2013}.

\subsection{Virtual elements on polyhedral cells}

Let $\P$ denote a generic three-dimensional element with boundary
$\partial\P$ and $\F\in\partial\P$ a generic polygonal face of $\P$.
We first introduce the elemental boundary space
\begin{align}
  \Bsh(\partial\P)
  := \Big\{
  \vsh\in\CS{0}(\partial\P):\,
  \restrict{\vsh}{\F}\in\Vsh(\F)
  \quad\forall\F\in\partial\P
  \Big\}.
  \label{eq:Bsh:def}
\end{align}
The functions in $\Bsh(\partial\P)$ are continuous linear polynomials
across the face edge and their restriction to a given face $\F$ belong
to the virtual element space $\Vsh(\F)$ defined
in~\eqref{eq:VEM:space:2D:enhanced}.
We use~\eqref{eq:Bsh:def} to introduce the preliminary virtual element
space on the polytopal element $\P$
\begin{align}
  \Vsht(\P)
  := \Big\{
  \vsh\in\HS{1}(\P):\,
  \restrict{\vsh}{\partial\P}\in\Bsh(\partial\P),\,
  \Delta\vsh\in\PS{1}(\P)
  \Big\}.
  \label{eq:VEM:space:3D:extended}
\end{align}
As for the two-dimensional case, the virtual element functions
$\vsh\in\Vsht(\P)$ are uniquely characterized by their vertex values,
which are again returned by the same linear bounded functionals
$\lambda_{\ell}(\vsh)$ introduced in the previous section.
We call these vertex values the \emph{degrees of freedom} of the
method.
The unisolvence of the set $\{\lambda_{\ell}(\vsh)\}$ in $\Vsh(\P)$ is
proved by using the same arguments as
in~\cite{Ahmad-Alsaedi-Brezzi-Marini-Russo:2013}.
Next, we introduce the elliptic projection operator
$\PinP{}:\Vsht(\P)\to\PS{1}(\P)$, which is such that
\begin{align}
  \int_{\P}\nabla\big(\vsh-\PinP{1}\vsh\big)\cdot\nabla\qs &= 0 \quad\forall\qs\in\PS{1}(\F),\label{eq:PinP:def:A}\\[0.5em]
  \int_{\partial\P}\big(\vsh-\PinP{1}\vsh\big)               &= 0.\label{eq:PinP:def:B}
\end{align}
As for the two-dimensional case, the elliptic operator $\PinP{1}\vsh$
of a virtual element function $\vsh\in\Vsht(\P)$ only depends on its
vertex values, i.e., the values $\lambda_{\ell}(\vsh)$.
Then, we define the local virtual element space
\begin{align}
  \Vsh(\P)
  := \bigg\{
  \vsh\in\Vsht(\P):\,\int_{\P}\vsh\qs=\int_{\P}\big(\PinP{1}\vsh\big)\qs
  \quad\forall\qs\in\PS{1}(\P)
  \bigg\}.
\end{align}
A major property of definition~\eqref{eq:VEM:space:2D:enhanced} is
that the linear polynomials over $\P$ are a subspace of $\Vsh(\P)$;
formally, we can write that
$\PS{1}(\P)\subseteq\Vsh(\F)\subseteq\Vsht(\P)$.

\medskip
It is worth noting that the orthogonal projection operator
$\PizP{\ell}:\Vsh(\P)\to\PS{\ell}(\P)$, $\ell=0,1$, is also computable
and, in particular, $\PizP{1}$ coincides with $\PinP{1}$,
see~\cite{Ahmad-Alsaedi-Brezzi-Marini-Russo:2013}.
The orthogonal projection $\PizP{\ell}\vsh$ of the virtual element
function $\vsh$ is the linear polynomial solving the variational
problem:
\begin{align}
  \int_{\P}\Big(\PizP{\ell}\vsh-\vsh\Big)\qs\dV=0
  \quad\forall\qs\in\PS{\ell}(\P).
\end{align}
Building on top of this definition, we can introduce the global
piecewise linear and constant orthogonal projections onto the
piecewise discontinuous polynomial space $\PS{\ell}(\Th)$ that is such
that $\restrict{\Piz{\ell}\vsh}{\P}=\PizP{\ell}(\restrict{\vsh}{\P})$
for all $\P\in\Th$.

\medskip
Finally, the global virtual element space $\Vsh\subset\HS{1}(\Omega)$
is defined by ``gluing'' all the elemental virtual element spaces in a
conforming way
\begin{align}
  \Vsh
  := \Big\{
  \vsh\in\HS{1}(\Omega):\,
  \restrict{\vsh}{\P}\in\Vsh(\P)
  \quad\forall\P\in\Th
  \Big\}.
  \label{eq:VEM:space:3D:enhanced}
\end{align}
The degrees of freedom that uniquely characterize the virtual element
function are still the vertex values $\lambda_{\ell}(\vsh)$, and their
unisolvence in the global space $\Vsh$ is a consequence of their
unisolvence in each elemental space.

\subsection{Virtual element approximation of the Poisson equation}

The virtual element approximation of~\eqref{eq:varform} is the
variational problem that reads as
\begin{align}
  \mbox{\textit{Find $\ush\in\Vshg$ such that}}\quad
  \ash(\ush,\vsh) = \bil{\fsh}{\vsh}
  \quad\forall\vsh\in\Vshzr,
  \label{eq:varform:VEM}
\end{align}
where
\begin{itemize}
\item $\Vshg:=\big\{\vsh\in\Vsh:\,\restrict{\vsh}{\Gamma}=\gsI\big\}$ and
  $\Vshzr:=\big\{\vsh\in\Vsh:\,\restrict{\vsh}{\Gamma}=0\big\}$;
\item $\gsI$ is a suitable approximation (trace interpolation) of the
  boundary function used in~\eqref{eq:Poisson:b};
\item the bilinear form $\ash:\Vsh\times\Vsh\to\REAL$ is the virtual
  element approximation of the left-hand side of
  Eq.~\eqref{eq:Poisson:a};
\item the linear functional $\bil{\fsh}{\cdot}:\Vsh\to\REAL$ is the
  virtual element approximation of the right-hand side of
  Eq.~\eqref{eq:Poisson:a} using $\fsh$, which is an element of the
  dual space $(\Vsh)^*$.
\end{itemize}
The well-posedness of this problem follows from the coercivity and
continuity of the bilinear form $\ash(\cdot,\cdot)$ and the
boundedness on the linear functional $\bil{\fsh}{\cdot}$, and can be
proved by an application of the Lax-Milgram theorem,
see~\cite[Theorem~2.7.7]{Brenner-Scott:2008}.
These properties follows from the construction that is briefly
reviewed below.

The bilinear form $\ash(\cdot,\cdot)$ is written as the sum of local
terms
\begin{align}
  \ash(\ush,\vsh) = \sum_{\P\in\Th}\ashP(\ush,\vsh),
\end{align}
where each term $\ashP:\Vsh(\P)\times\Vsh(\P)\to\REAL$ is a symmetric
bilinear form over the elemental space $\Vsh(\P)$.
We set
\begin{align}
  \ashP(\ush,\vsh) 
  = \asP(\PinP{1}\ush,\PinP{1}\vsh) 
  + \hP\SP(\ush-\PinP{1}\ush,\vsh-\PinP{1}\vsh),
  \label{eq:poly:ah:def}
\end{align}
and the stabilization form $\SP:\Vsh\times\Vsh\to\REAL$ is a
symmetric, positive definite, bilinear form such that
\begin{align}
  \sigma_*\asP(\vsh,\vsh)\leq\SP(\vsh,\vsh)\leq\sigma^*\asP(\vsh,\vsh)
  \quad\forall\vsh\in\Vsh\textrm{~with~}\PinP{1}\vsh=0,
  \label{eq:poly:S:stability:property}
\end{align}
for two positive constants $\sigma_*$, $\sigma^*$ that are independent
of $\hh$ (and $\P$).
Condition~\eqref{eq:poly:S:stability:property} implies that the
resulting bilinear form $\ashP(\cdot,\cdot)$ in~\eqref{eq:poly:ah:def}
has the two major properties of \emph{linear consistency} and
\emph{stability}.
These two properties are exploited in the convergence analysis of the
method,
see~\cite{BeiraodaVeiga-Brezzi-Cangiani-Manzini-Marini-Russo:2013,Ahmad-Alsaedi-Brezzi-Marini-Russo:2013}.
In all our numerical test case, the solver implements the so
called~\emph{dofi-dofi} stabilization,
see~\cite{Dassi-Mascotto:2018,Mascotto:2018}.

\medskip
To approximate the right-hand side of \eqref{eq:varform}, we first set
$\restrict{\fsh}{\P}=\PizP{1}\fs$ and then we consider the elemental
decomposition
\begin{align}
  \bil{\fsh}{\vsh} = \sum_{\P\in\Th}\int_{\P}\big(\PizP{1}\fs)\vsh.
  \label{eq:vem:rhs}
\end{align}
Estimates of such approximation are
found in~\cite[Section~5.8]{Ahmad-Alsaedi-Brezzi-Marini-Russo:2013}.

\subsection{Convergence results}
\label{subsec:vem:convergence}
\medskip
For the sake of reference, we report below the convergence result in
the $\LS{2}$-norm and the energy norm for the numerical approximation
using the virtual element space~\eqref{eq:VEM:space:3D:enhanced}.
This result follows from the general convergence theorem that is
proved in~\cite[Theorem~1]{Ahmad-Alsaedi-Brezzi-Marini-Russo:2013}.
Let $\us\in\HS{2}(\Omega)$ be the solution to the variational
problem~\eqref{eq:varform} on a convex domain $\Omega$ with
$\fs\in\HS{1}(\Omega)$.
Let $\ush\in\Vsh$ be the solution of the virtual element method
\eqref{eq:varform:VEM} on every mesh of a mesh family $\calT=\{\Th\}$
satisfying a suitable set of mesh geometrical assumptions.
Then, a strictly positive constant $\Cs$ exists such that
\begin{itemize}
\item the $\HONE$-error estimate holds:
  \begin{align}
    \label{eq:source:problem:H1:error:bound}
    \norm{\us-\ush }{1,\Omega}\leq
    \Cs\hh\left(
    \norm{\us}{2,\Omega} 
    + \snorm{\fs}{1,\Omega}
    \right);
  \end{align}
\item the $\LTWO$-error estimate holds:
  \begin{align}
    \label{eq:source:problem:L2:error:bound}
    \norm{\us-\ush}{0,\Omega}\leq 
    \Cs\hh^{2}\left(
    \norm{\us}{2,\Omega}
    + \snorm{\fs}{1,\Omega}
    \right).
  \end{align}
\end{itemize}
Finally, we note that the approximate solution $\ush$ is not
explicitly known inside the elements.
Consequently, in the numerical experiments of
Section~\ref{sec:results}, we approximate the error in the
$\LTWO$-norm as follows:
\begin{align}
  \label{eq:error:H1}
  \norm{\us-\ush}{0,\Omega}\approx\norm{\us-\Piz{1}\ush}{0,\Omega},
\end{align}
where $\Piz{1}\ush$ is the global $\LTWO$-orthogonal projection of the
virtual element approximation $\ush$ to $\us$.
On its turn, we approximate the error in the energy norm as follows:
\begin{align}
  \label{eq:error:L2}
  \snorm{ \us-\ush
  }{1,\Omega}\approx\norm{\nabla\us-\Piz{0}\nabla\ush}{0,\Omega},
\end{align}
where $\Piz{0}$ is extended component-wisely to the vector fields.

\section{Analysis of the quality of a polyhedral mesh}
\label{sec:quality}
As mentioned in Section~\ref{subsec:vem:convergence}, there are
several theoretical results behind the VEM that depend on particular
geometrical (or \textit{regularity}) assumptions.
These assumptions ensure that 
in the refinement process each element of any mesh in the mesh family, is sufficiently regular,
and guarantee the VEM convergence and optimal estimates of the
approximation error with respect to different norms.
In this section we overview the most common geometrical assumptions
found in the literature and provide an indicator of the violation of
these assumptions, which depends uniquely on the geometry of the mesh
elements.
As these assumptions have been shown to be quite restrictive
\cite{sorgente2021role}, we try to measure \textit{how much} a mesh
satisfies them instead of simply discarding all the non-compliant meshes.

\subsection{Geometrical assumptions}
\label{subsec:assumptions}
A complete analysis of the geometrical assumptions typically required
in literature for the VEM in two dimensions can be found in
\cite{sorgente2021vem}.
We report here the main results from that paper, as they will be the
basis on which we build their three-dimensional counterparts.

\medskip
Besides some small variants, we can isolate four main conditions
common to all the considered works relative to the 2D case.
They are defined for a single mesh $\Th$, but the conditions contained
in them are required to hold independently of $h$.
Therefore, when considering a mesh family $\calT=\{\Th\}_h$, these
assumptions have to be verified simultaneously by every $\Th \in
\calT$.
We use the superscript 2D to indicate the fact
that they are relative to the two-dimensional VEM formulation. 
\begin{itemize}
\item[\GAs{1}]\textit{ $\exists\,\rho\in(0,1)$ s.t. every
  polygon $\P\in\Th$ with diameter $\hP$ is star-shaped with respect
  to a disc with radius
    \[\rP\ge\rho\hP.\]}
\item[\GAs{2}]\textit{ $\exists\,\rho\in(0,1)$ s.t. for every
  polygon $\P\in\Th$, the length $\hE$ of every edge $\E\in\bP$
  satisfies
    \[\hE\ge\rho\hP.\]}
\item[\GAs{3}]\textit{ $\exists\,\Ns\in\mathbb{N}$ s.t. the
  number of edges of every polygon $\P\in\Th$ is (uniformly) bounded
  by $\Ns$.\\}
\item[\GAs{4}]\textit{ $\exists\,\delta>0$ s.t. for every
  polygon $\P\in\Th$ the one-dim. mesh $\calI_\P$ representing its
  boundary can be subdivided into a finite number of disjoint
  sub-meshes $\calI_\P^1, \ldots, \calI_\P^N$ (each one containing
  possibly more than one edge of $\P$), and for each $\calI_\P^i$ it
  holds that
  \smalleq{\frac{\max_{\E\in\calI_\P^i}\hE}{\min_{\E\in\calI_\P^i}\hE}
    \le \delta.}}
\end{itemize}
Authors typically require their meshes to satisfy either
\GAs{1} and \GAs{2} or \GAs{1} and
\GAs{2} or \GAs{1} and \GAs{3}, also depending
to the type of stabilization term adopted.

\medskip
The assumptions for the regularity of polyhedral meshes are straightforwardly derived from the 2D ones.
According to 
\cite{Ahmad-Alsaedi-Brezzi-Marini-Russo:2013,%
BeiraodaVeiga-Brezzi-Marini-Russo:2014,%
BeiraodaVeiga-Dassi-Russo:2017,%
Brenner:2017:SEV,%
brenner2018virtual},
we consider the following requirements.
\begin{itemize}
\item[\GAv{1}]\textit{ $\exists\,\rho\in(0,1)$ s.t. every
  polyhedron $\P\in\Th$ is star-shaped with respect to a disc with
  radius
  \[\rP\ge\rho\hP,\]
  and every face $\F\in\partial\P$ is star-shaped with respect to a
  disc with radius
  \[\rF\ge\rho\hF,\]}
\item[\GAv{2}]\textit{ $\exists\,\rho\in(0,1)$ s.t. for every
  polyhedron $\P\in\Th$, the length $\hE$ of every edge $\E$ of every
  face $\F$ satisfies
  \[\hE\ge\rho\hF\ge\rho^2\hP.\]}
\item[\GAv{3}]\textit{ $\exists\,\Ns\in\mathbb{N}$ s.t. the number
  of edges and faces of every polyhedron $\P\in\Th$ is (uniformly)
  bounded by $\Ns$.}
\end{itemize}
In most of the considered papers only \GAv{1} and \GAv{2} are considered, while more recent works like \cite{brenner2018virtual} replace \GAv{2} with the more general \GAv{3}.
Indeed, \GAv{3} can be derived from \GAv{1} and
\GAv{2}, as these latters imply the existence of an integer number
$\Ns$ such that every polyhedron has fewer than $\Ns$ faces and every
face has fewer than $\Ns$ edges \cite[Remark~11]{Ahmad-Alsaedi-Brezzi-Marini-Russo:2013}.
It is worth noting that here we are not
extending~\cite[Assumption~\textbf{G4}]{sorgente2021role} to the
three-dimensional case.
To the best of our knowledge such extension has not yet been
considered in the analysis of the 3D VEM, and this topic is beyond the
scope of our work.

As observed in \cite[Remark~10]{Ahmad-Alsaedi-Brezzi-Marini-Russo:2013},
such assumptions allow us to use very general mesh families.
However, as for the two-dimensional case, they can be more restrictive
than necessary to ensure the proper convergence of the method.
This fact is also something that we want to investigate in our work.

\subsection{Quality indicators}
\label{subsec:indicator}
In the first paper \cite{sorgente2021role}, starting from each geometrical assumption \GAs{i}, \textbf{i~}$=1,\ldots, 4$, we derived a scalar function $\rhos{i}: \{\P \subset \Omega_h\} \to [0,1]$ defined element-wise, which measures how well a polygon $\P\in \Omega_h$ meets the requirements of \GAs{i} from 0 ($\P$ does not respect \GAs{i}) to 1 ($\P$ fully respects \GAs{i}).
\begin{align}
\rhos{1}(\P) &= \frac{k(\P)}{\mP} = 
    \begin{cases} 
    1 & \mbox{if $\P$ is convex} \\
    \in (0,1) & \mbox{if $\P$ is concave and star-shaped} \\
    0 & \mbox{if $\P$ is not star-shaped} \\
    \end{cases}\\
\rhos{2}(\P) &= \frac{\min(\sqrt{\mP}, \ \min_{\E\in\partial\P}\hE)}{\max(\sqrt{\mP}, \ \hP)} \\
\rhos{3}(\P) &= \frac{3}{\#\left\{\E\in\partial\P\right\}} =
    \begin{cases} 
    1 & \mbox{if $\P$ is a triangle} \\
    \in (0,1) & \mbox{otherwise} \\
    \end{cases}\\
\rhos{4}(\P) &= \min_i \frac{\min_{\E\in \calI_\P^i} \hE}{\max_{\E\in\calI_\P^i} \hE}
\end{align}
Combining together $\rhos{1}, \rhos{2}$, $\rhos{3}$ and $\rhos{4}$, we defined a global function $\rhos{}:\{\Th\}_h\to [0,1]$ which measures the overall quality of a mesh $\Omega_h$:
\begin{small}
\begin{equation}
\rhos{}(\Omega_h) = \sqrt{\frac{1}{\#\left\{\P\in\Omega_h \right\}} \sum_{\P\in\Omega} \frac{\rhos{1}(\P)\rhos{2}(\P) + \rhos{1}(\P)\rhos{3}(\P) + \rhos{1}(\P)\rhos{4}(\P)}{3}}.
\label{eq:indicator2d}
\end{equation}
\end{small}
We have $\rhos{}(\Omega_h)=1$ if and only if $\Omega_h$ is made only of equilateral triangles, $\rhos{}(\Omega_h)=0$ if and only if $\Omega_h$ is made only of non star-shaped polygons, and $0<\rhos{}(\Omega_h)<1$ otherwise.

\medskip
Similarly, for the 3D case, from each geometrical assumption \GAv{i}, \textbf{i}$=1,2,3$, we need to derive a scalar function $\rhov{i}$ which measures how well a polyhedron $\P\in \Omega_h$ meets the requirements of \GAv{i}.
We do not consider $\rhos{4}$, since assumption \GAs{4} has not been extended to the 3D scenario.
We measure the quality of the interior of $\P$ defining a new volumetric operator, and we also include a measure of the quality of its faces $\partial\P$ using the old 2D indicators.
Then we collect together the volumetric and the boundary parts into a single function.

\begin{align}
\rhov{1}(\P) &= \frac{k(\P)}{\mP}
    \prod_{\F\in\partial\P}\rhos{1}(\F) \\
    &= \begin{cases} 
    1 & \mbox{if $\P$ and all its faces are convex} \\
    \in (0,1) & \mbox{if $\P$ and all its faces are concave and star-shaped} \\
    0 & \mbox{if $\P$ or one of its faces are not star-shaped} \\
    \end{cases}\notag\\
\rhov{2}(\P) &= \frac{1}{2}\sqbra{
    \frac{\min(\sqrt[3]{\mP},\,\min_{\F\in\partial\P}\hF)}{\max(\sqrt[3]{\mP},\,\hP)} +
    \frac{1}{\#\cubra{\F\in\partial\P}} \sum_{\F\in\partial\P}\rhos{2}(\F)} \\
\rhov{3}(\P) &= \frac{1}{2}\sqbra{
    \frac{4}{\#\cubra{\F\in\partial\P}} +
    \frac{1}{\#\cubra{\F\in\partial\P}} \sum_{\F\in\partial\P}\rhos{3}(\F)}\\
    &=\begin{cases} 
    1 & \mbox{if $\P$ is a tetrahedron} \\
    \in (0,1) & \mbox{otherwise} \\
    \end{cases}\notag
\end{align}
The operator $k(\cdot)$ in $\rhov{1}$ measures the volume of the \textit{kernel} of a polyhedron $\P$, defined as the set of points in $\P$ from which the whole polyhedron is visible, computed with the algorithm proposed in \cite{SorgenteKernel}.
The volumetric and the boundary components (the kernel of each face) are multiplied so that, even if only one of them is zero (if a single face is not star-shaped), the whole product vanishes.
The function $\rhov{2}$ is an average of the volumetric constant $\rho$ from \GAv{2}, expressed trough the ratio $\hF/\hP$, and all the boundary constants represented by the 2D indicators $\rhos{2}$.
Function $\rhov{3}$ measures the number of edges and faces of a polyhedron, penalizing elements with numerous edges and faces as required by \GAv{3}.

We can now define a global function $\rhov{}:\{\Th\}_h\to [0,1]$ which measures the overall quality of a polyhedral mesh $\Omega_h$.
The formula for combining $\rhov{1}$, $\rhov{2}$ and $\rhov{3}$ is derived straightforwardly from \eqref{eq:indicator2d}:
%
%
\begin{equation}
\rhov{}(\Omega_h) = \sqrt{\frac{1}{\#\left\{ \P \in \Omega_h \right\}} \ \sum_{\P\in\Omega} \frac{\rhov{1}(\P)\rhov{2}(\P) + \rhov{1}(\P)\rhov{3}(\P)}{2}}.
\label{eq:indicator3d}
\end{equation}
We have $\rhov{}(\Omega_h)=1$ if and only if $\Omega_h$ is made only of equilateral tetrahedrons, $\rhov{}(\Omega_h)=0$ if and only if $\Omega_h$ is made only of non star-shaped polyhedrons (or with non star-shaped faces), and $0<\rhov{}(\Omega_h)<1$ otherwise.
All indicators $\rhov{1}, \rhov{2}$ and $\rhov{3}$, and consequently $\rhov{}$, only depend on the geometrical properties of the mesh elements; therefore their values can be computed before applying the VEM, or any other numerical scheme.

\medskip
As already observed in \cite{sorgente2021role}, this approach is easily upgradeable to future developments: whenever new assumptions on the features of a mesh should appear in the literature (for example a valid extension of \GAs{4}), one simply needs to introduce in this framework a new function $\rhov{i}$ that measures the violation of the new assumption and inserts it into the formulation of the general indicator $\rhov{}$ in equation~(\ref{eq:indicator3d}).
For the rest of the work we will be dealing only with volumetric meshes, therefore we will indicate $\rhov{}$ simply as $\varrho$ and \GAv{i} as \textbf{Gi} for \textbf{i}$=1,2,3$, with a little abuse of notation.

\section{Generation of the datasets}
\label{sec:datasets}

In this section, we define a number of mesh ``datasets'' on which we will test our indicator.
We call a \emph{dataset} a collection $\calD:=\{\Omega_n\}_{n=0,\ldots,N}$ of meshes $\Omega_n$ covering the domain $\Omega=(0,1)^3$ such that:
\begin{itemize}
    \item the mesh $\Omega_{n+1}$ has smaller mesh size than $\Omega_{n}$ for every $n=0,\ldots,N-1$;
    \item the meshes $\Omega_n$ are built with the same technique and therefore contain similar polyhedra organized in similar configurations.
\end{itemize}
Note that each mesh $\Omega_n$ is uniquely identified (in the dataset it belongs to) by its size as $\Th$, therefore we can consider a dataset $\calD$ as a subset of a mesh family: $\calD=\{\Th\}_{h\in\calH'}\subset\calT$ where $\calH'$ is a finite subset of $\calH$.

\medskip
Each dataset is characterized by a \textit{sampling strategy} for generating a number of points inside the unit cube $\Omega$, and a \textit{meshing technique} to connect them.
Algorithms for generating all the samplings strategies and the meshing techniques are based on \textit{cinolib} \cite{livesu2019cinolib}; all the generated datasets are publicly available at:\\
\centerline{\url{https://github.com/TommasoSorgente/vem-indicator-3D-dataset}}

\subsection{Sampling strategies}
\label{subsec:dataset:sampling}
\begin{figure}[ptb]
\centering
\begin{tabular}{ccc}
\includegraphics[width=.2\linewidth]{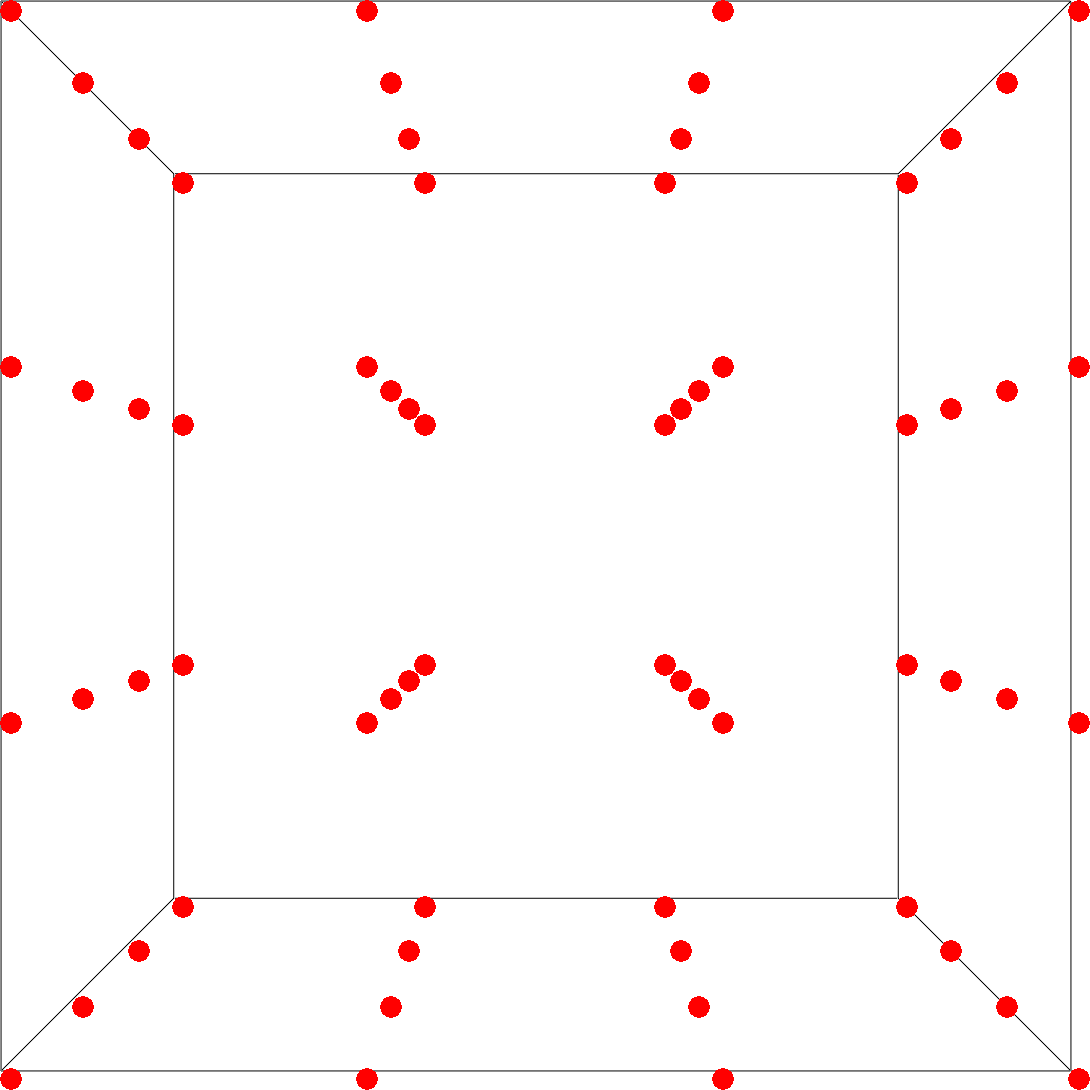} &
\includegraphics[width=.2\linewidth]{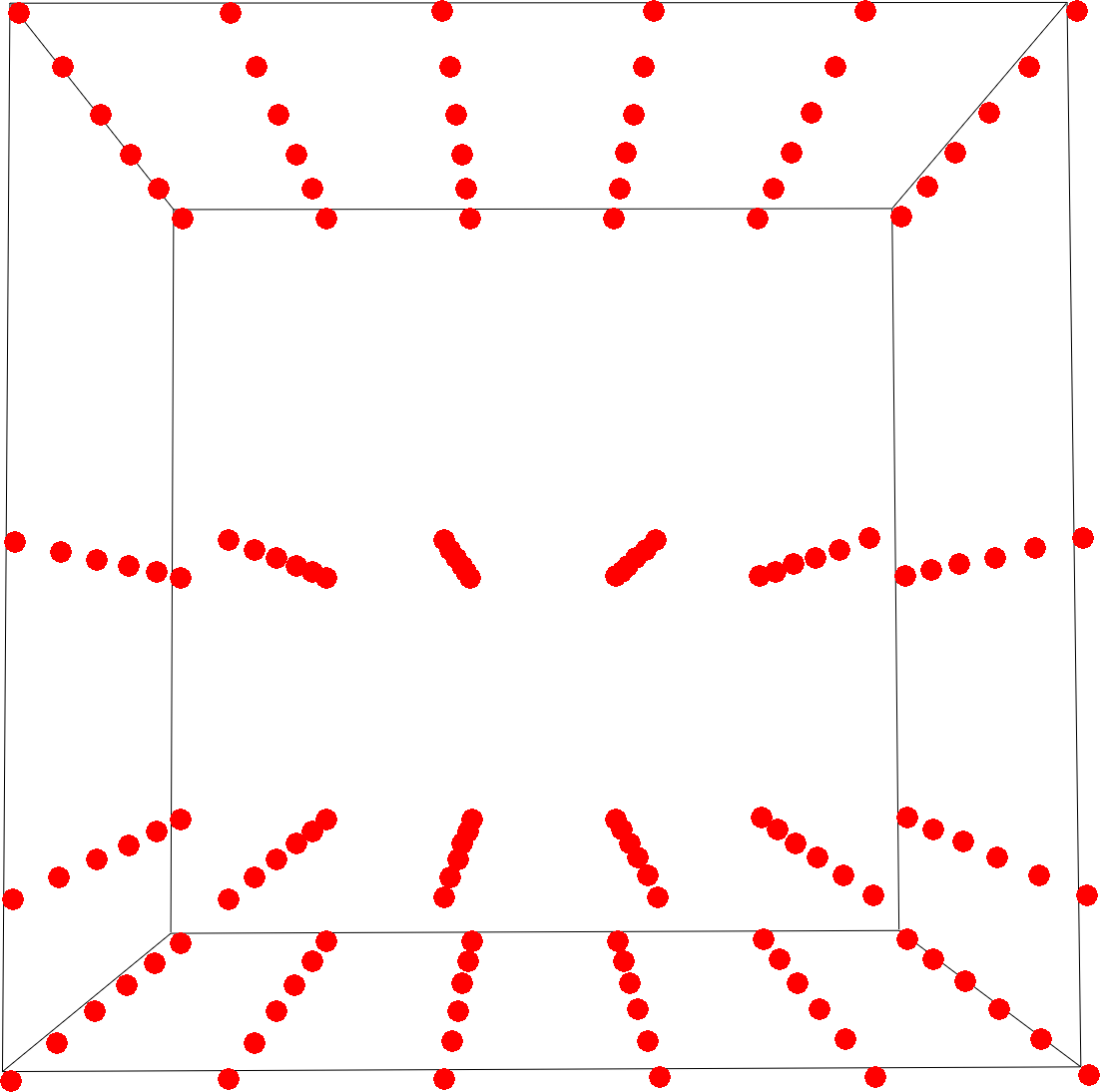} &
\includegraphics[width=.2\linewidth]{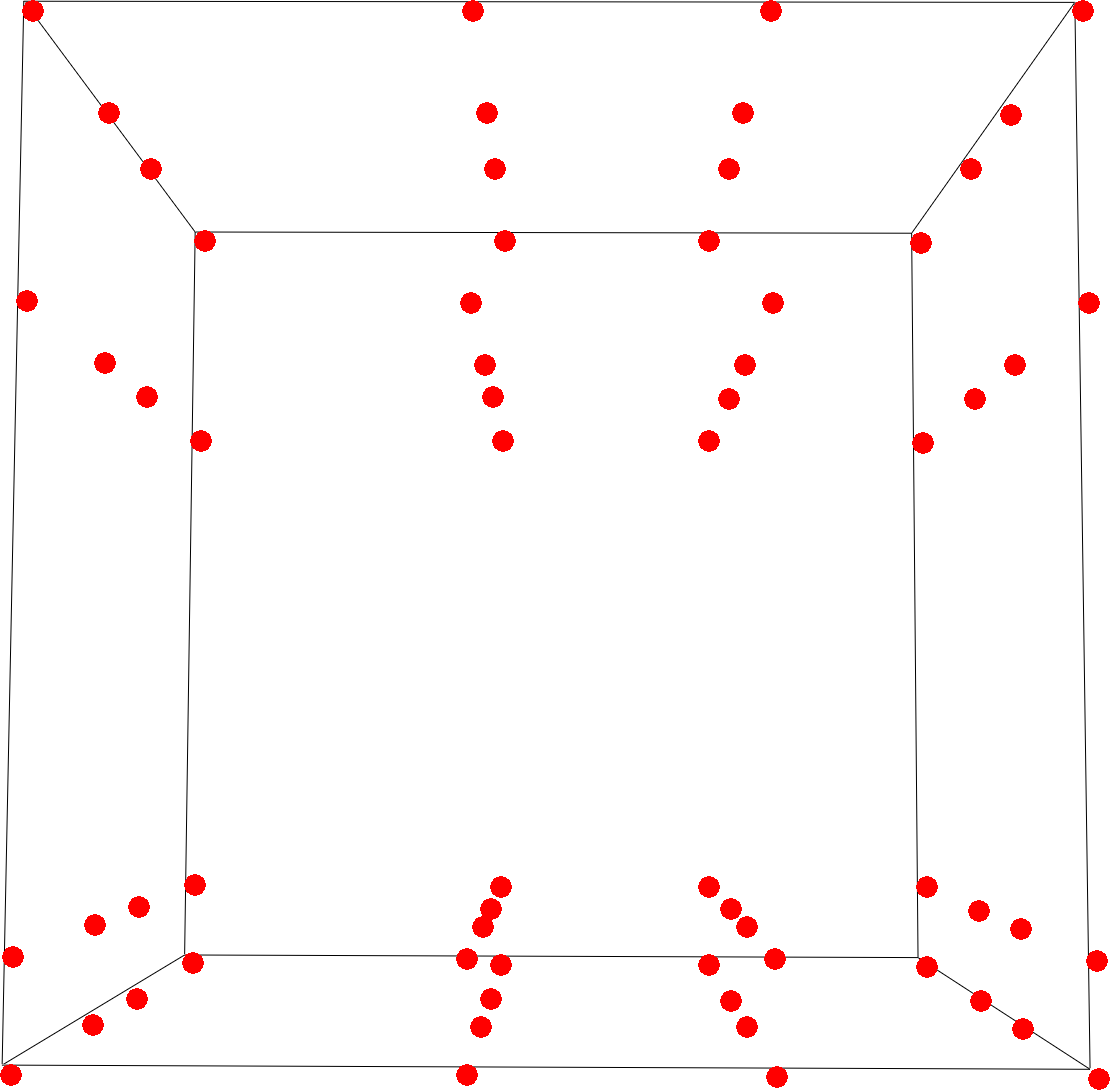} \\
{\footnotesize uniform sampling} & {\footnotesize anisotropic sampling} & {\footnotesize parallel sampling} \\
\vspace{0.2cm}\\
\includegraphics[width=.2\linewidth]{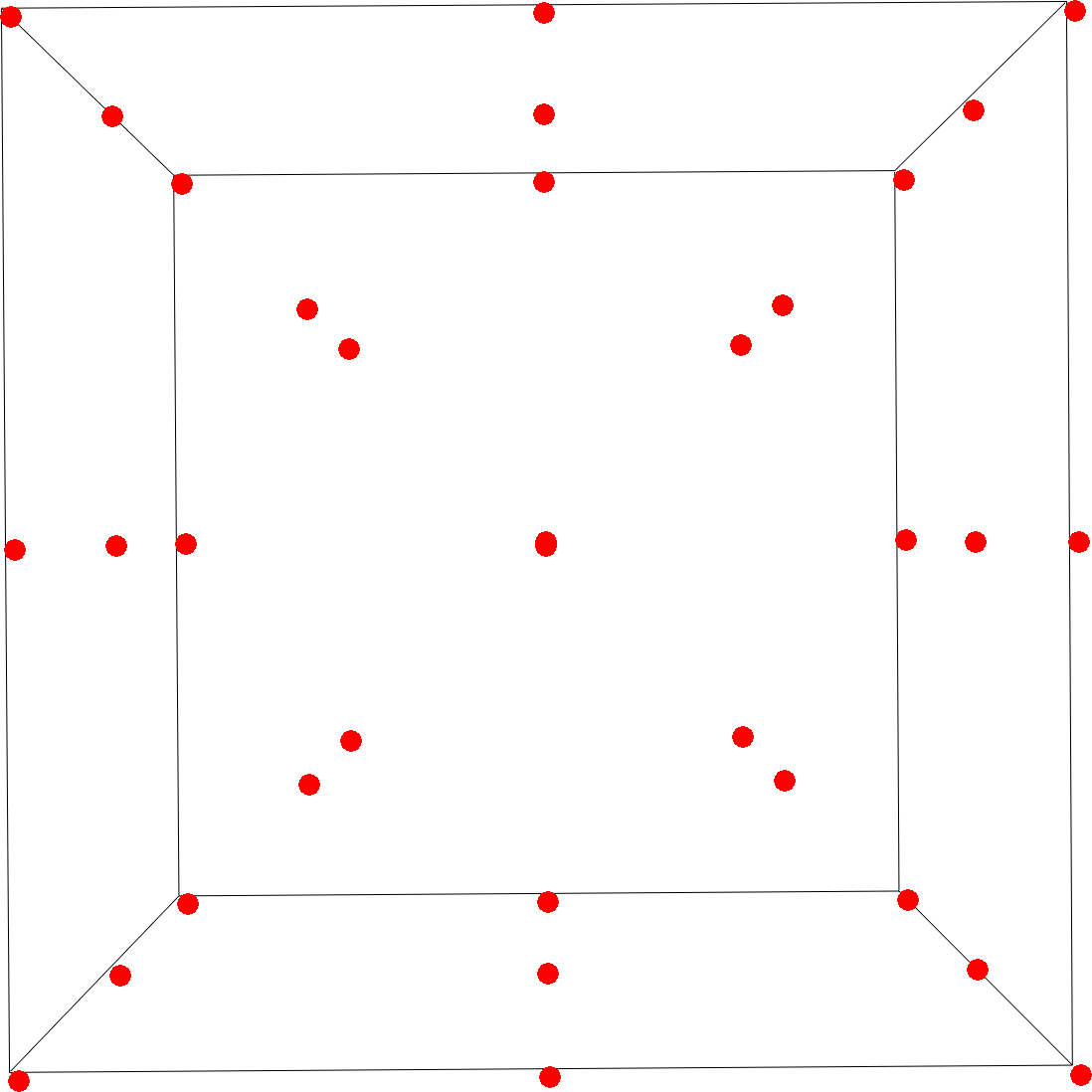} &
\includegraphics[width=.2\linewidth]{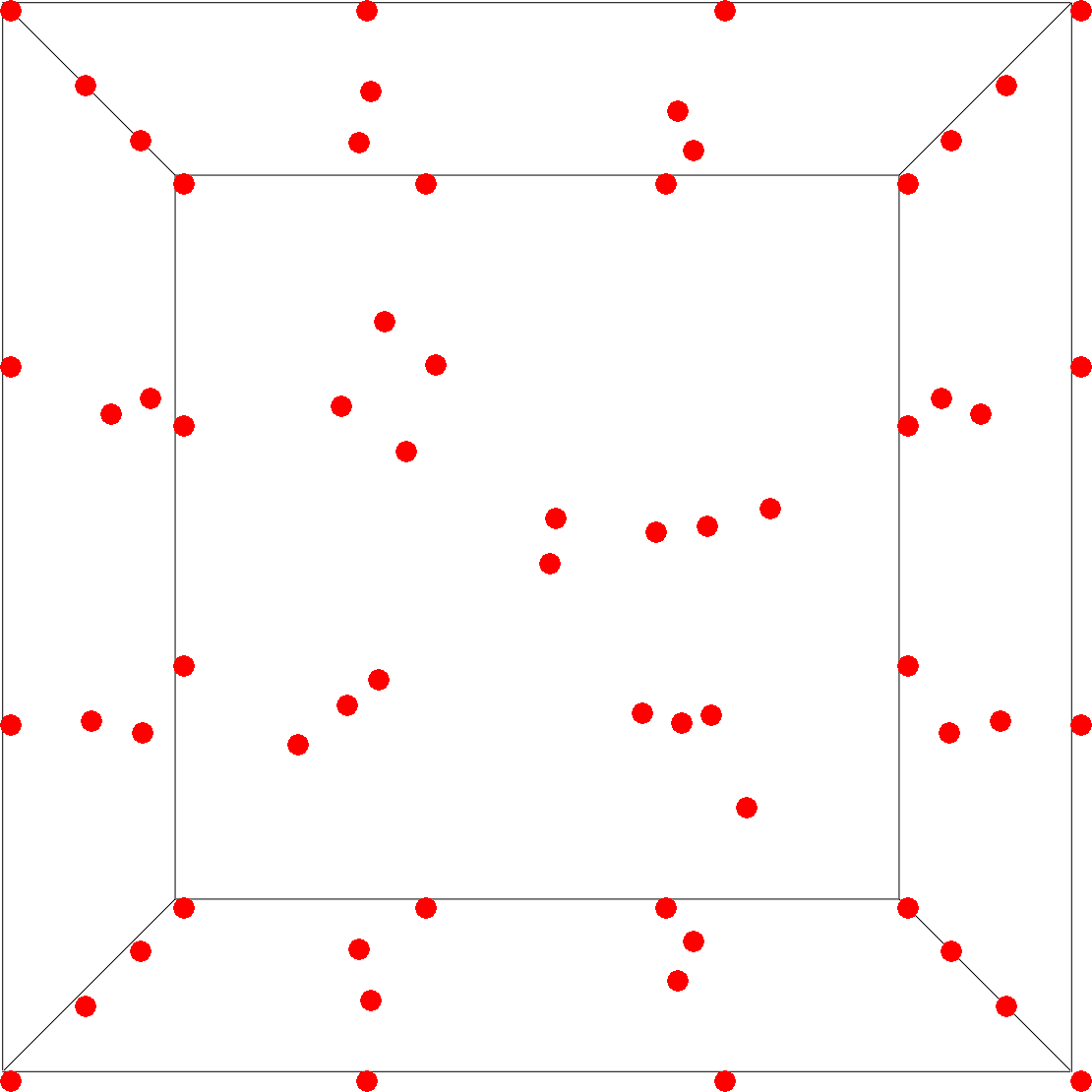} & 
\includegraphics[width=.2\linewidth]{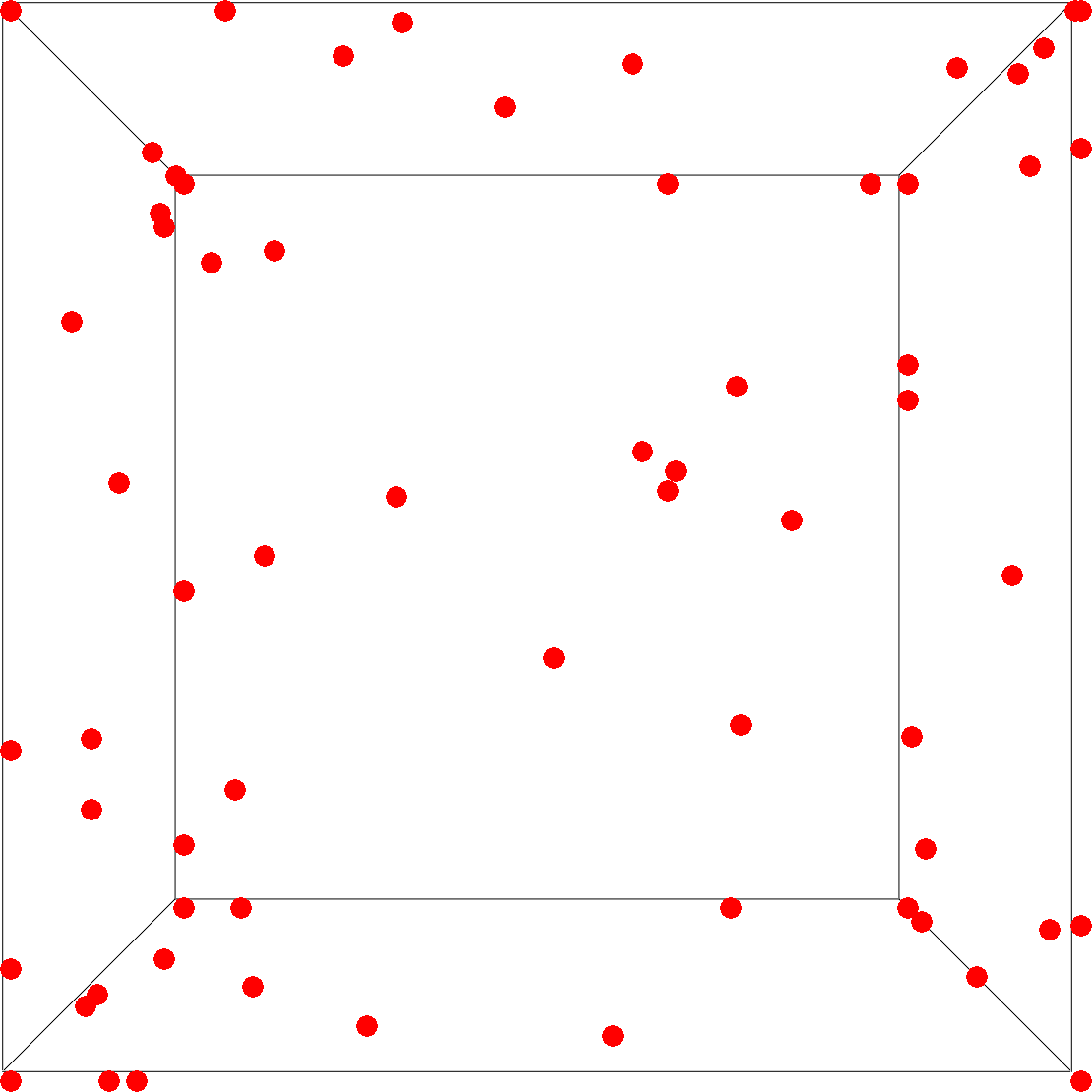}\\
{\footnotesize Body Centered Lattice} & {\footnotesize Poisson sampling} & {\footnotesize random sampling}
\end{tabular}
\caption{Summary of the sampling strategies.}
\label{fig:points}
\end{figure}

We defined six different sampling strategies, summarized in Figure~\ref{fig:points}.
In order to create multiple refinements for each dataset, each sampling takes an integer parameter $t$ in input, which determines the number of points to be generated, and consequently the size of the induced mesh.

The sampling strategies are the followings:
\begin{itemize}
    \item{\it Uniform sampling} 
    the points are disposed along a uniform equispaced grid of size $1/t$;
    \item{\it Anisotropic sampling} 
    a regular grid in which the distance between the points is fixed at $1/t$ along two directions while it linearly increases ($1/t$, $2/t$, $3/t$, $\ldots$) along the third axis, leading to anisotropic configurations;
    \item{\it Parallel sampling} 
    this sampling is obtained from a uniform sampling with parameter $t$ by randomly moving all the points belonging to a certain plane $p_1$ to another plane $p_2$, parallel to $p_1$.
    In practice, we pick all the points sharing the same $x-$coordinate and randomly change $x$ by the same quantity; then, we repeat this operation for the $y$ and the $z-$coordinates;
    \item{\it Body Centered Lattice (BCL)} 
    a uniform equispaced grid of size $1/t$ is generated and one more point is added to the center of each cubic cell, using the implementation proposed in \url{https://github.com/csverma610/CrystalLattice}.
    This sampling produces equilateral tetrahedra when combined with tetrahedral meshing and truncated octahedra when combined with Voronoi meshing (see Section~\ref{subsec:dataset:meshing});
    \item{\it Poisson sampling} 
    the points are generated following the Poisson Disk Sampling algorithm \cite{bridson2007fast}.
    First, we apply the algorithm to the cube $(1/t, 1-1/t)^3$ and the square $(1/t, 1-1/t)^2$ with radius $1/t$, to generate points inside $\Omega$ and on its boundary.
    Then, in order to cover the domain more uniformly, we add an equispaced sampling with distance $1/t$ on each edge of $\Omega$;
    \item{\it Random sampling} 
    the points are randomly placed inside $\Omega$. 
    In order to guarantee a decent distribution, given the input parameter $t$ we generate $t$ points along each edge of $\Omega$, $t^2$ points on each face and $t^3$ points inside $\Omega$.
\end{itemize}
We point out that, for how the samplings are defined, the number of points generated with the same parameter $t$ in the different strategies may vary.

\subsection{Meshing techniques}
\label{subsec:dataset:meshing}
We create meshes by connecting a certain set of points with different techniques: we considered the three most common types of mesh connectivity found in the literature (tetrahedral, hexahedral and Voronoi) plus a generic polyhedral one.
As soon as sampling points are connected into a mesh, we call them \textit{vertices}.

\begin{itemize}
    \item{\it Tetrahedral meshing} 
    points are connected in tetrahedral elements with \textit{TetGen} library \cite{si2015tetgen}, enforcing two quality constraints on tetrahedra: a maximum radius-edge ratio bound and a minimum dihedral angle bound;
    \item{\it Voronoi meshing} 
    points are considered as centroids for the construction of a Voronoi lattice using library \textit{Voro++} \cite{rycroft2009voro}.
    In this case the sampled points will not appear as vertices in the final mesh, and the number of vertices does not depend uniquely on the number of points;
    \item{\it Hexahedral meshing} 
    points are connected to form hexahedral elements.
    For ease of implementation the mesh is still generated as a Voronoi lattice, but we make sure that points are placed in such a way that the final result is a pure hexahedral mesh;
    \item{\it Polyhedral meshing} 
    we start from a tetrahedral mesh and we aggregate $20\%$ of the elements to generate non-convex polyhedra.
    In order to avoid numerical problems, we select the elements with the greatest volume and aggregate them with the neighboring element sharing the widest face.
    Note that we only aggregate couple of elements and merge eventual coplanar faces, but in any case we remove vertices of the original tetrahedral mesh.
\end{itemize}
It would also be possible to generate polyhedral meshes starting from hexahedral or Voronoi ones, but we chose to limit our tests to aggregations of the tetrahedral datasets.

\subsection{Validation datasets}
\label{subsec:dataset:datasets}
\begin{figure}[ptb]
\centering
\begin{tabular}{ccc}
\includegraphics[width=.2\linewidth]{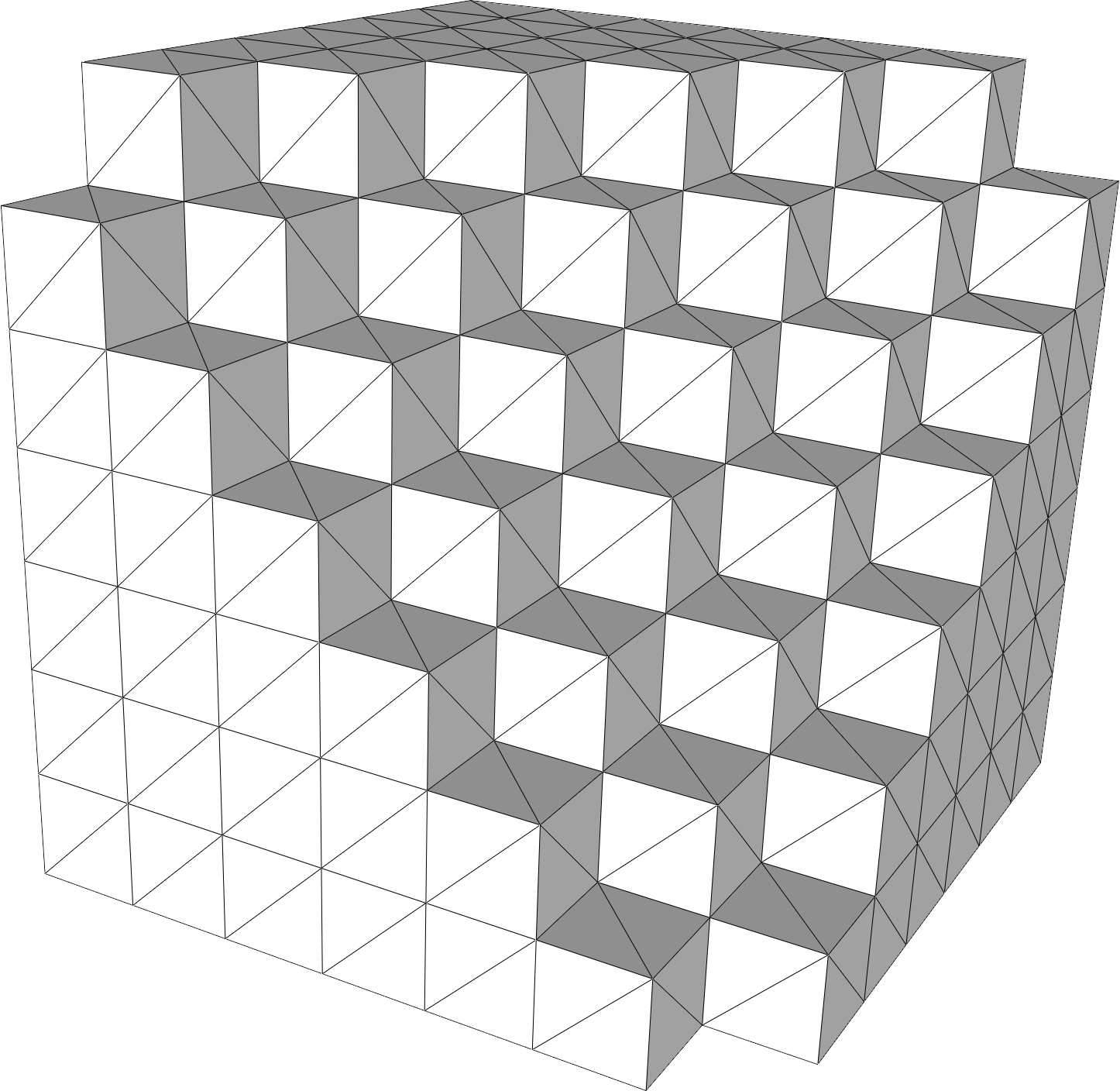} &
\includegraphics[width=.2\linewidth]{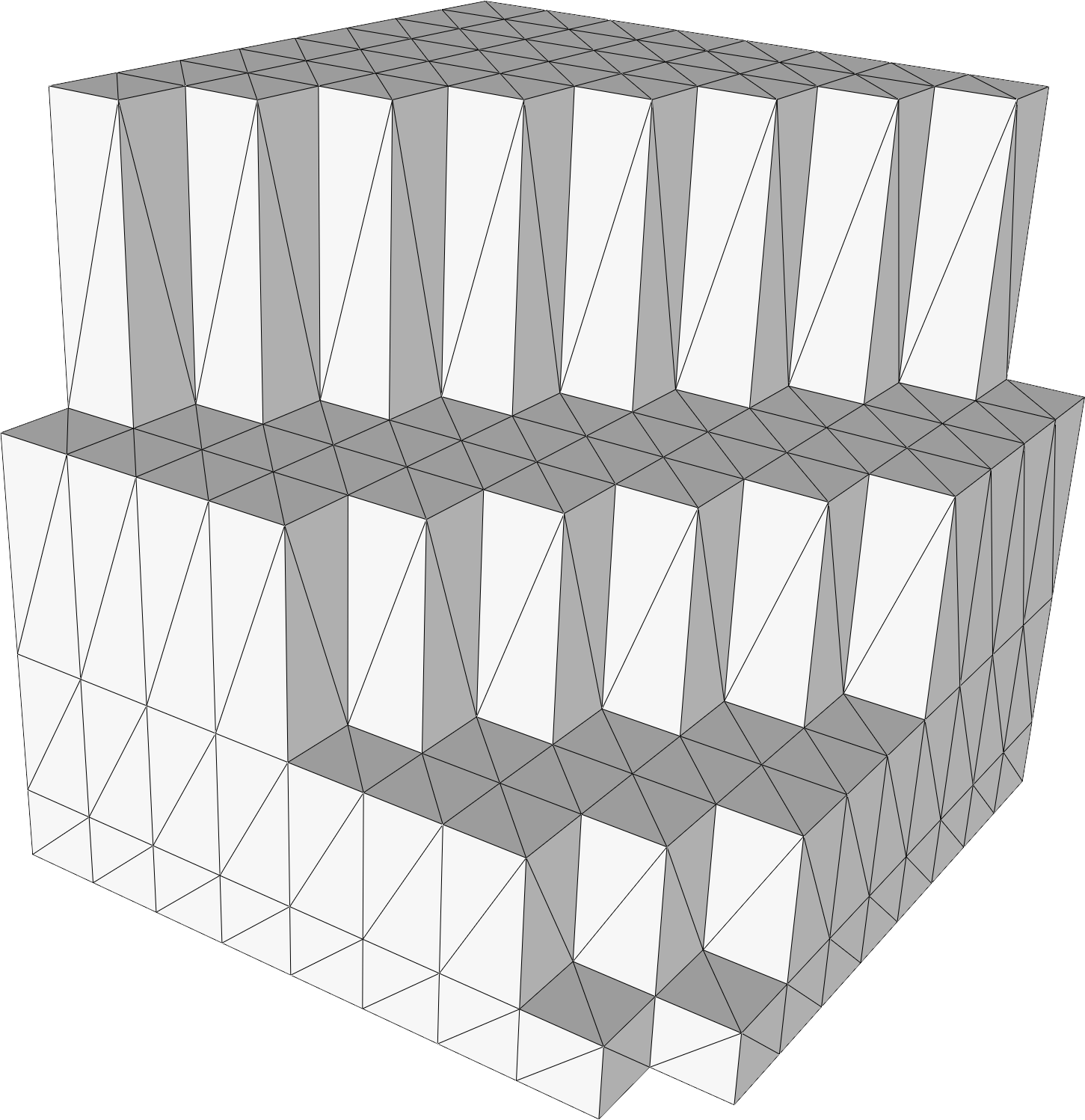} &
\includegraphics[width=.2\linewidth]{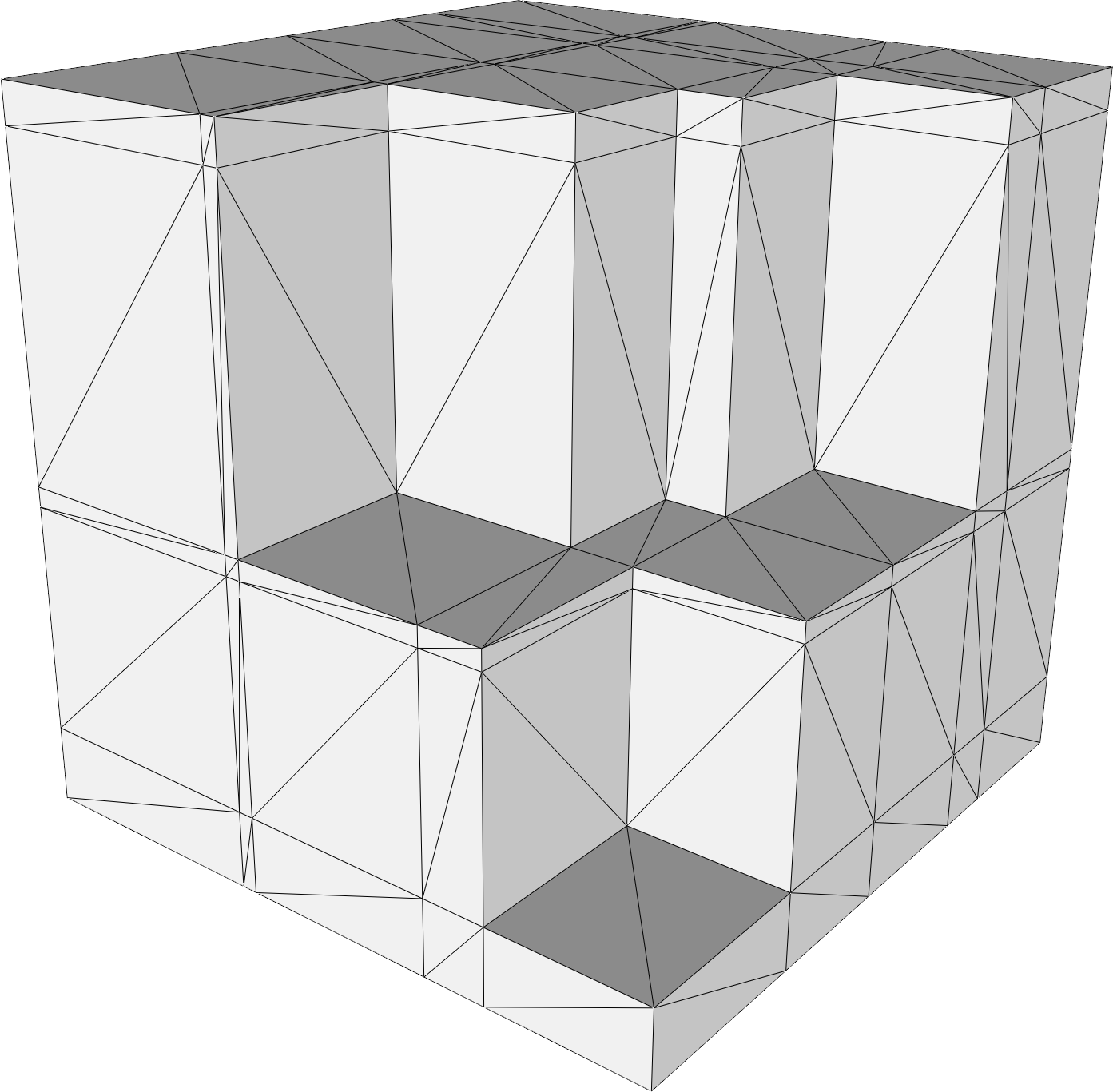}\\
$\Dtetuni$ & $\Dtetani$ & $\Dtetpar$ \\
\vspace{0.2cm}\\
\includegraphics[width=.2\linewidth]{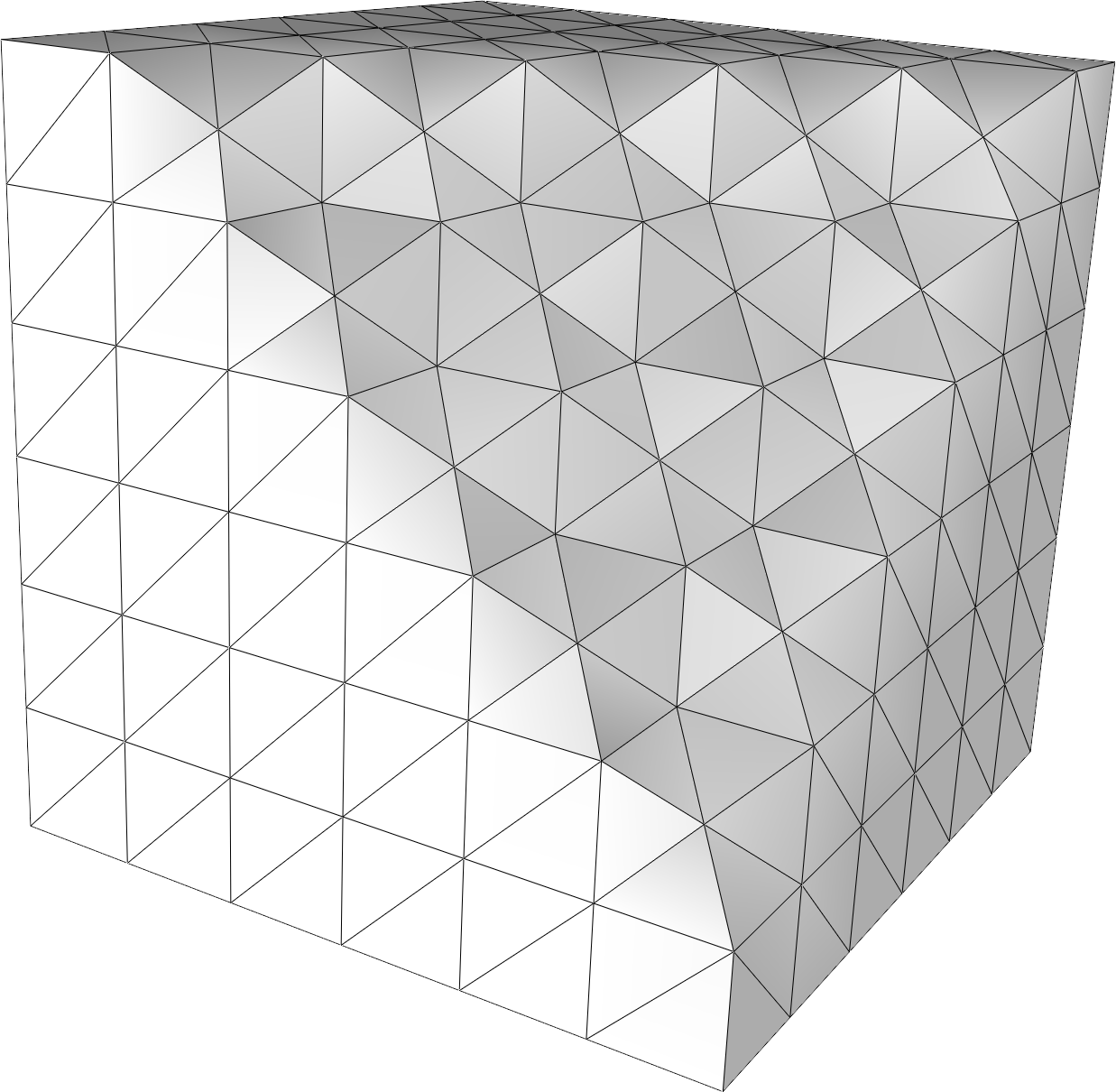} &
\includegraphics[width=.2\linewidth]{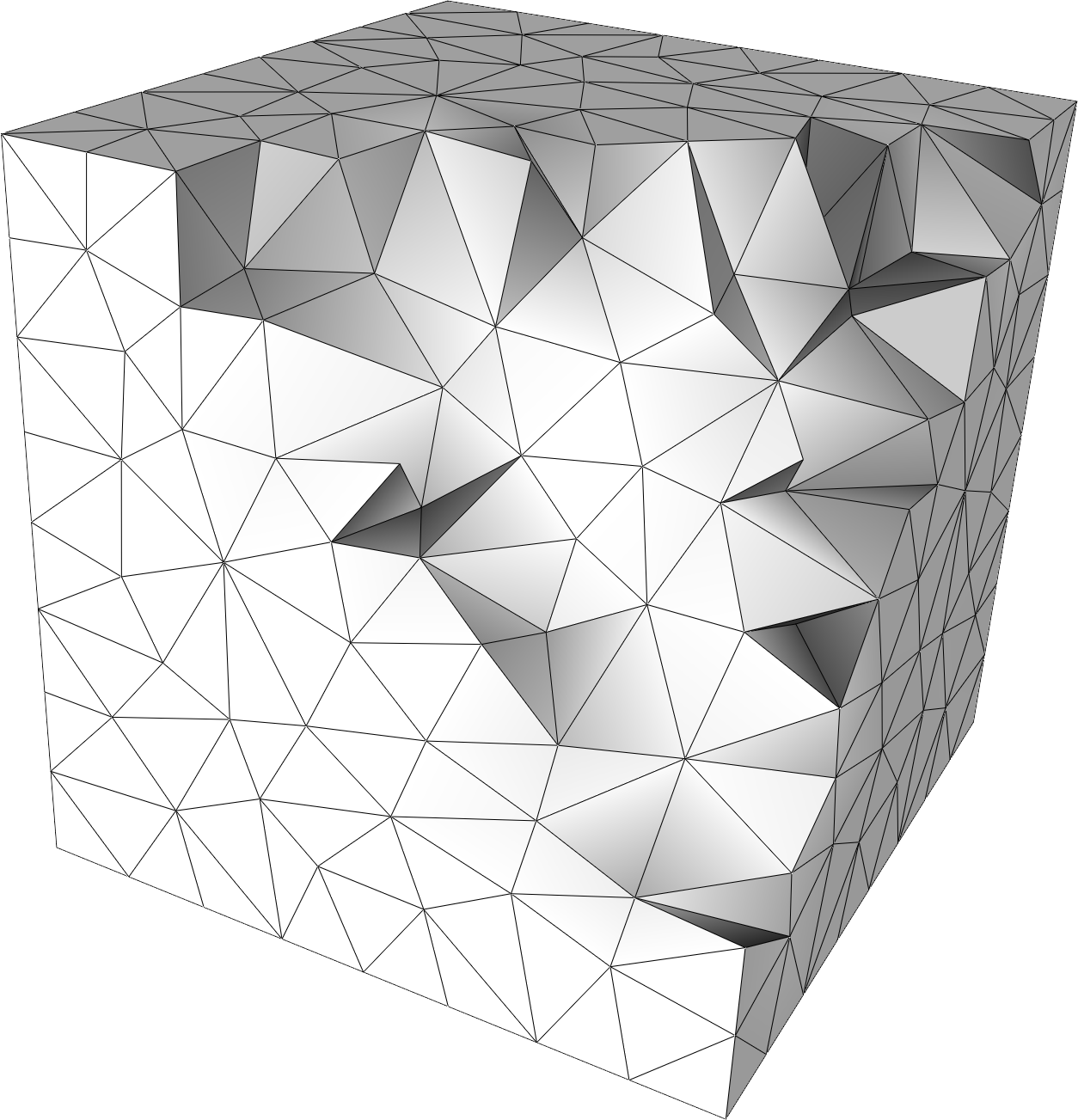} &
\includegraphics[width=.2\linewidth]{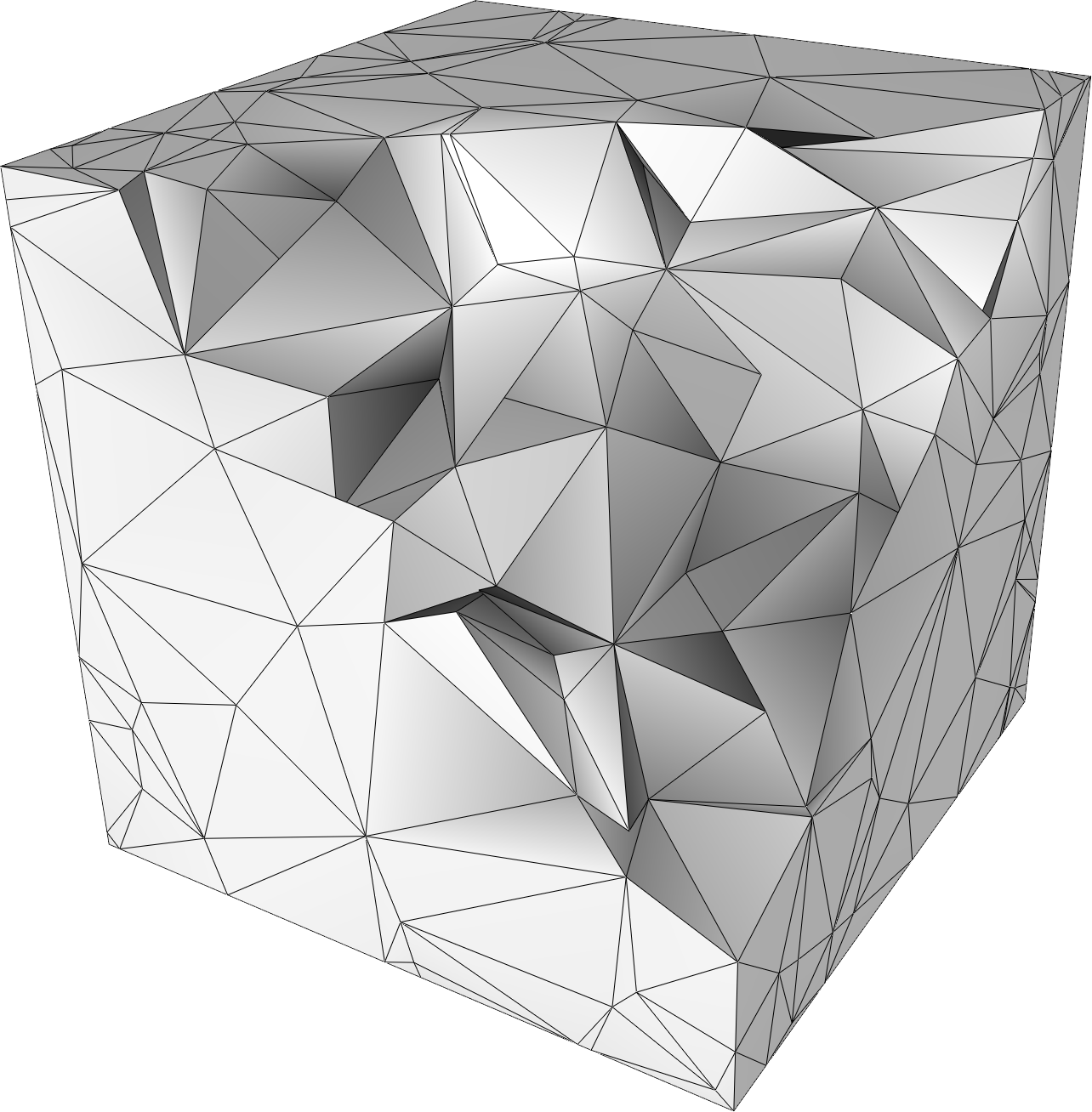} \\
$\Dtetbcl$ & $\Dtetpoisson$ & $\Dtetrandom$
\end{tabular}
\caption{Summary of the tetrahedral datasets}
\label{fig:tet-datasets}
\end{figure}

\begin{figure}[ptb]
\centering
\begin{tabular}{ccc}
\includegraphics[width=.2\linewidth]{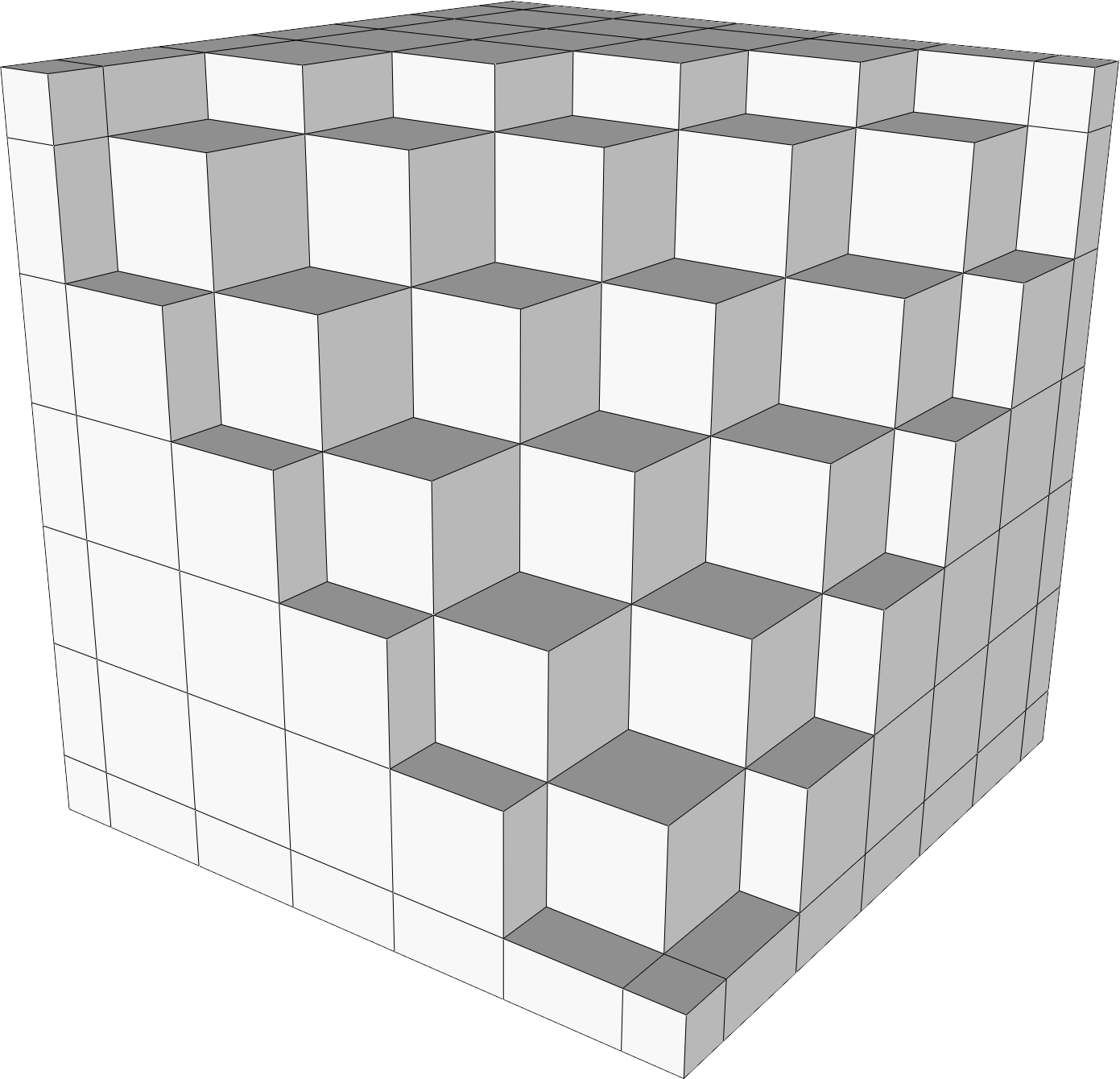} &
\includegraphics[width=.2\linewidth]{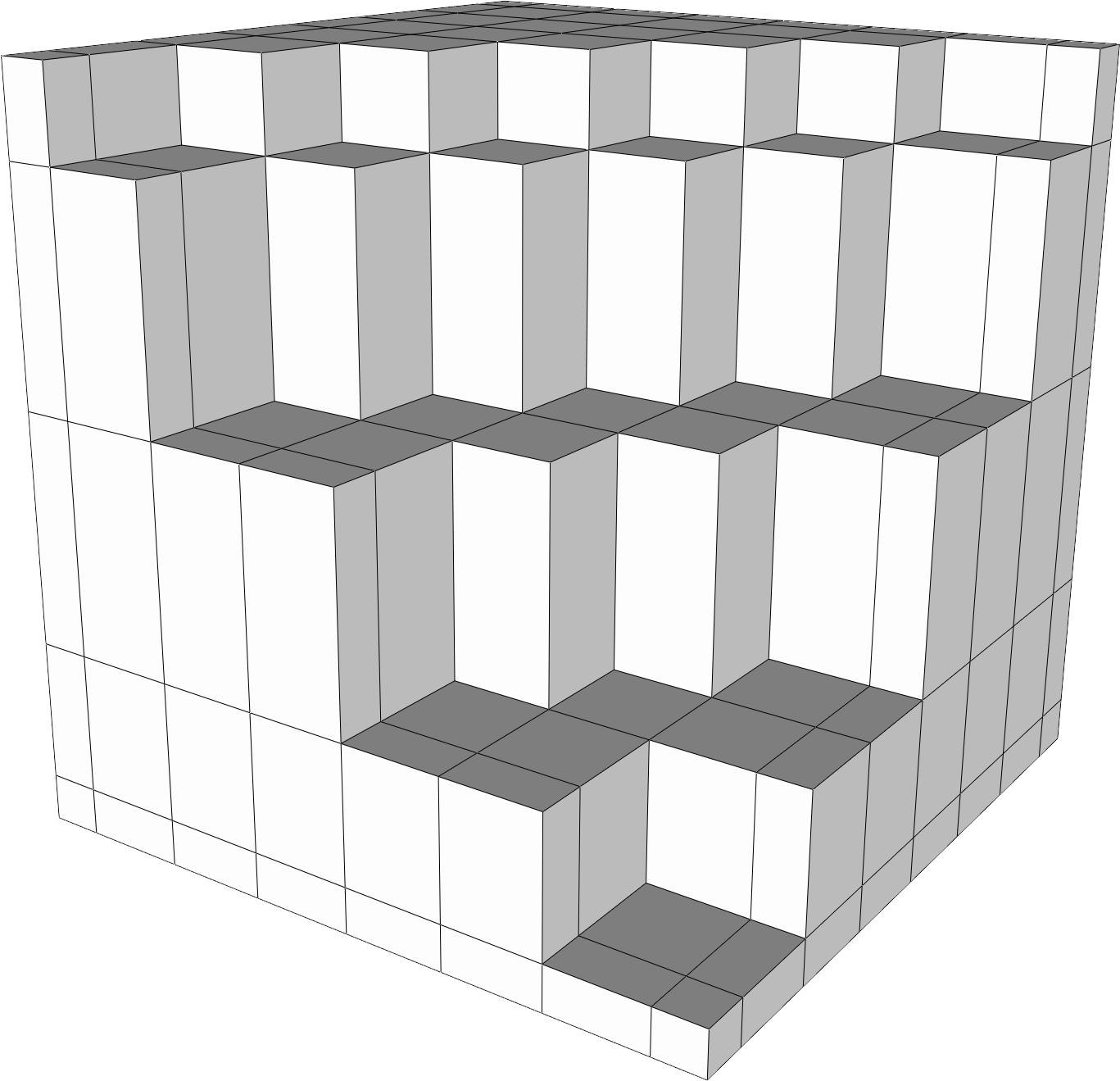} &
\includegraphics[width=.2\linewidth]{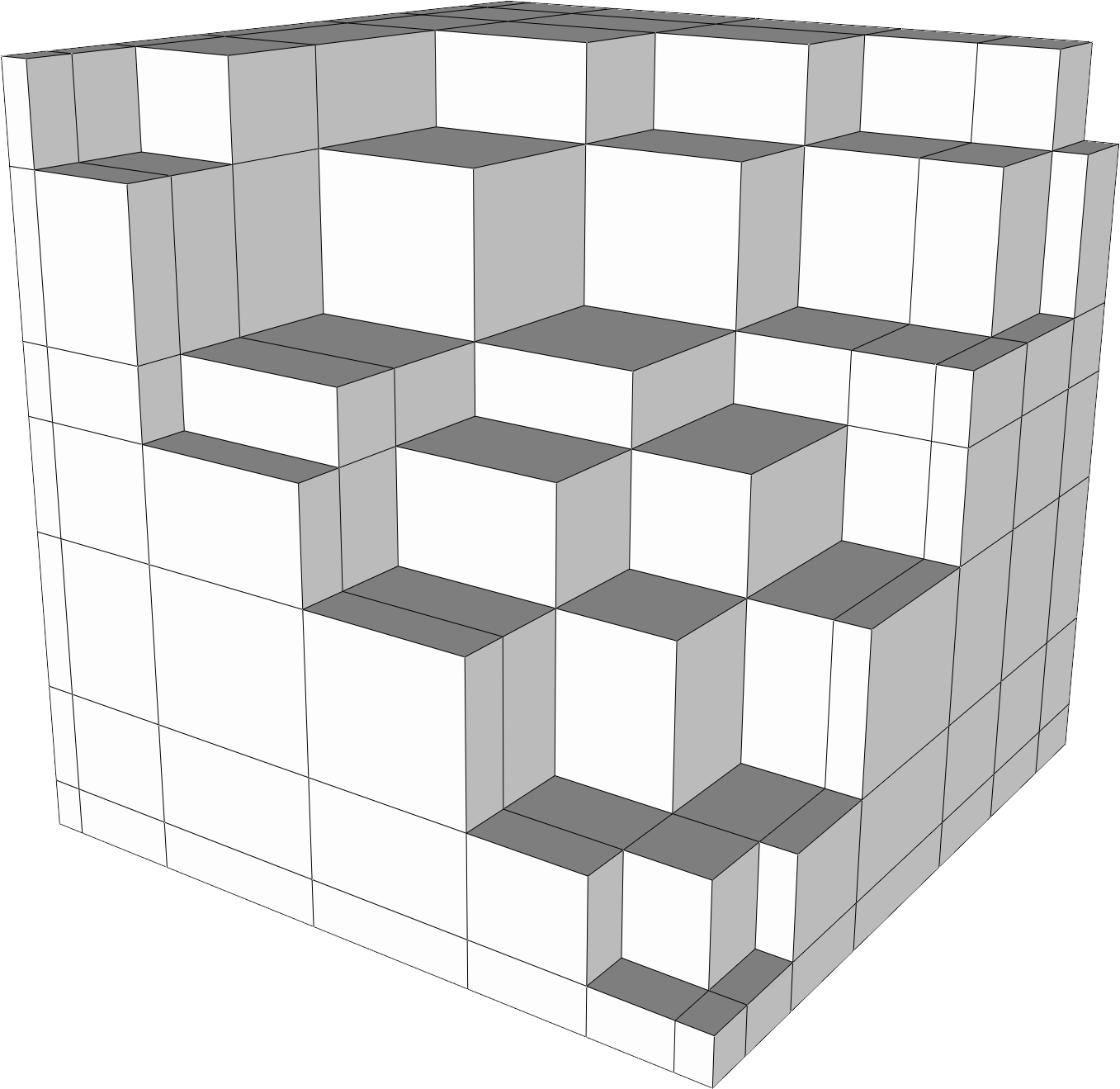}\\
$\Dhexuni$ & $\Dhexani$ & $\Dhexpar$ \\
\vspace{0.2cm}\\
\includegraphics[width=.2\linewidth]{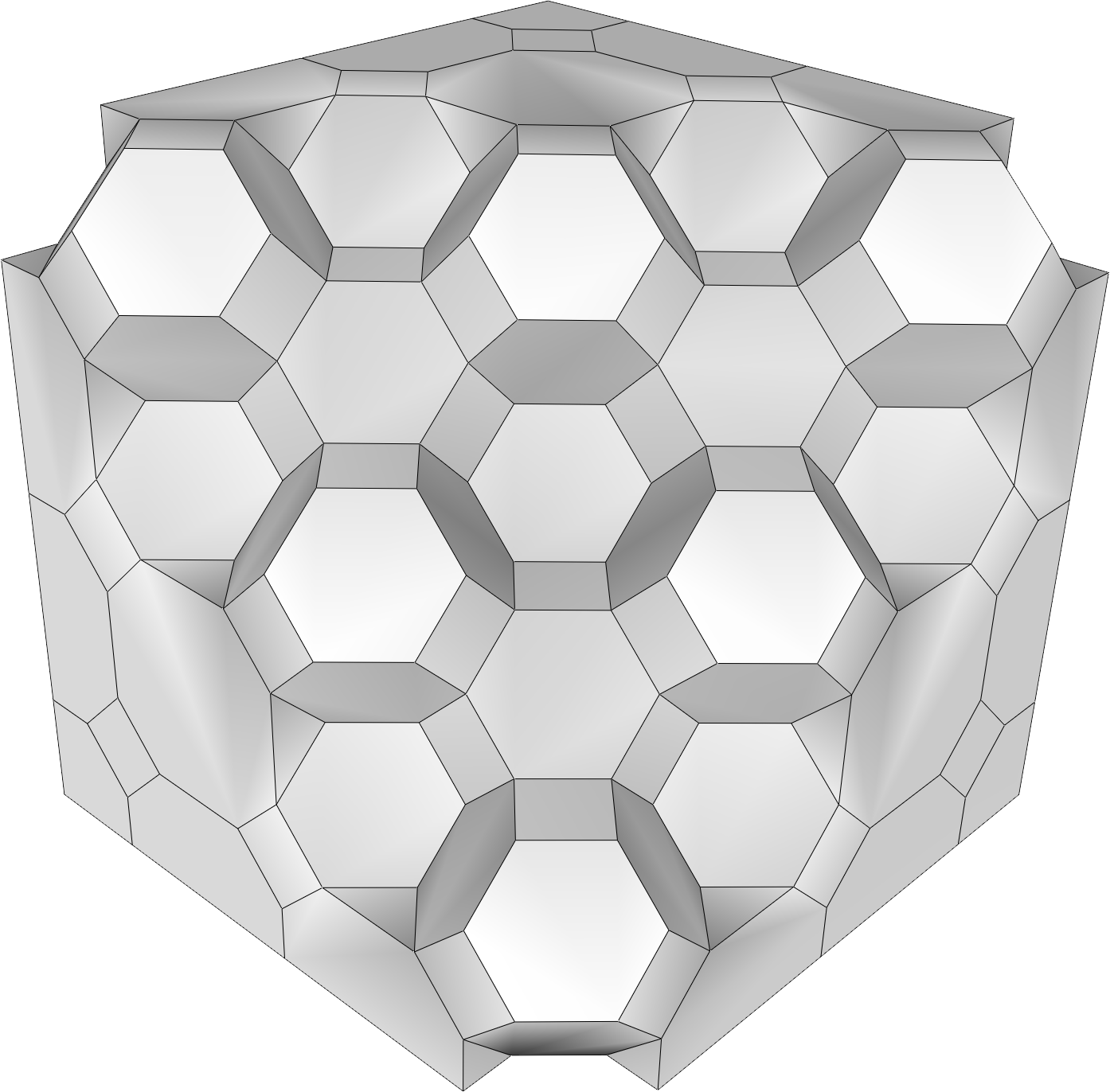} &
\includegraphics[width=.2\linewidth]{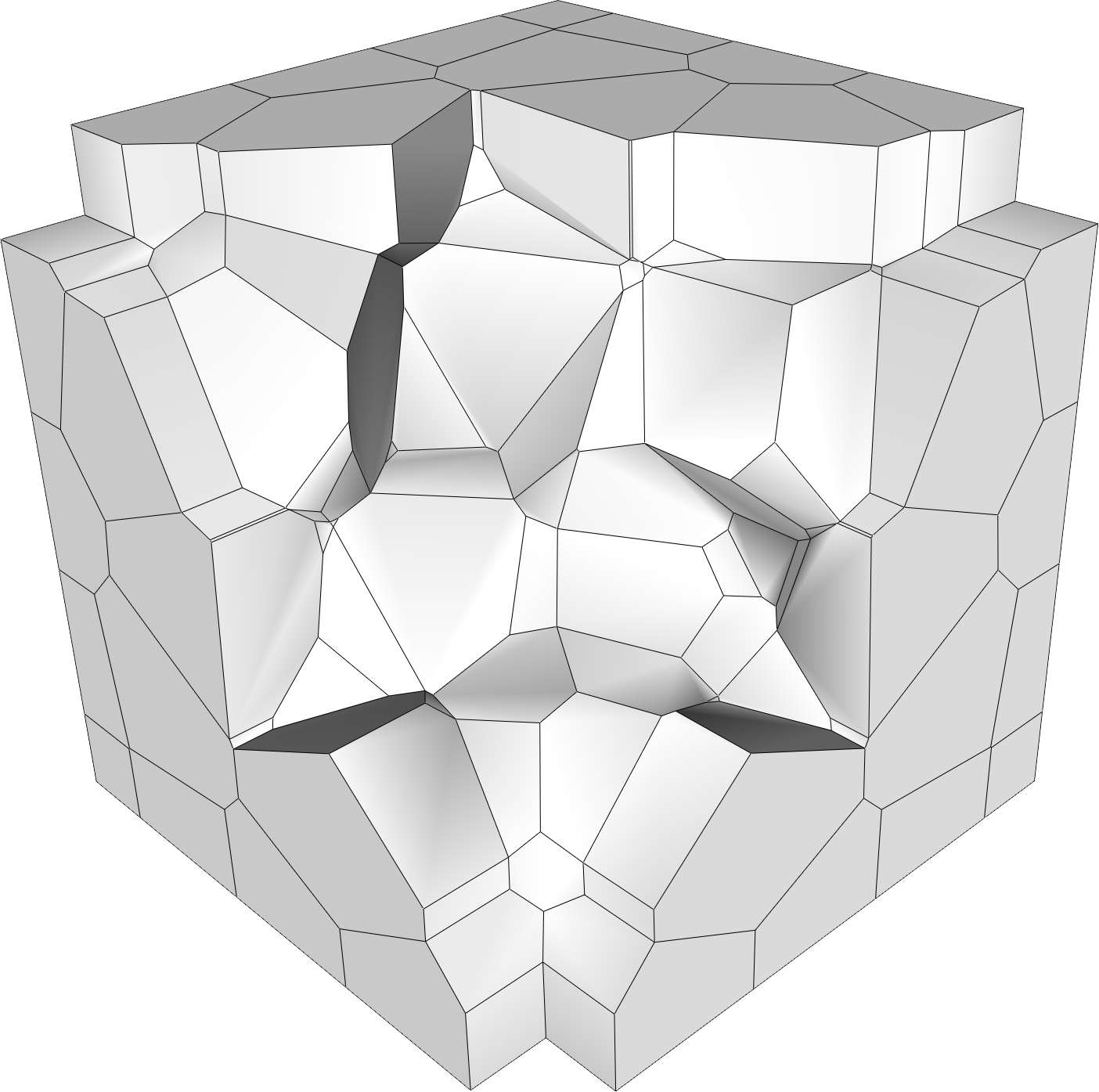} &
\includegraphics[width=.2\linewidth]{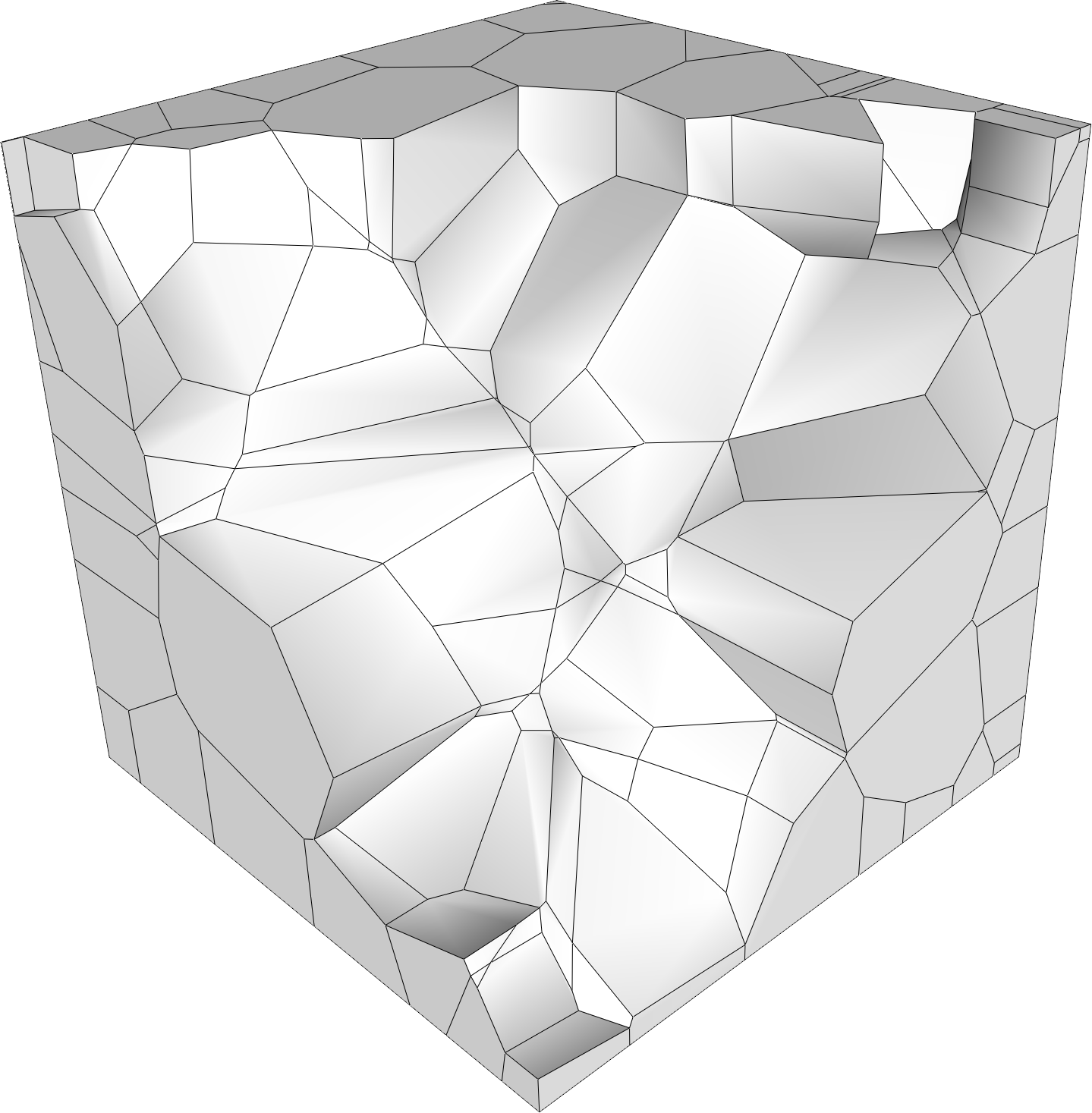} \\
$\Dvorobcl$ & $\Dvoropoisson$ & $\Dvororandom$
\end{tabular}
\caption{Summary of the hexahedral and Voronoi datasets}
\label{fig:hex-voro-datasets}
\end{figure}

Each combination of sampling strategy and meshing technique gives a dataset.
As not all combinations are possible or meaningful, we selected the most significant datasets for our purpose, presented in Figure~\ref{fig:tet-datasets} and Figure~\ref{fig:hex-voro-datasets}.

\medskip
For validating our experiments we generated datasets composed by five meshes each, with decreasing mesh size: they contain 60, 500, 4000, 32000 and 120000 vertices approximately.
For each $n$, the mesh $\Omega_{n+1}$ has about eight times more vertices than $\Omega_{n}$, except for the last mesh which only has four times more vertices than the previous, for computational reasons.
We determine the number of vertices by opportunely setting the sampling parameter $t$, and we have no constraints on the numbers of edges, faces or elements.
We label each dataset to indicate the meshing technique and the sampling strategy. For example, $\Dtetuni$ is the dataset that is built by combining the tetrahedral meshing and the uniform sampling.

\begin{itemize}
    \item{\it Tetrahedral datasets} 
    we combined the tetrahedral meshing with all the considered sampling methods, creating six tetrahedral datasets: $\Dtetuni$, $\Dtetani$, $\Dtetpar$, $\Dtetbcl$, $\Dtetpoisson$ and $\Dtetrandom$;
    \item{\it Hexahedral datasets} 
    among the six considered samplings, only the first three provide regular grids suitable for the generation of hexahedral elements.
    Therefore our hexahedral datasets are $\Dhexuni$, $\Dhexani$ and $\Dhexpar$;
    \item{\it Voronoi datasets} 
    the last three samplings instead, have been used to generate Voronoi datasets: $\Dvorobcl$, $\Dvoropoisson$ and $\Dvororandom$;
    \item{\it Polyhedral datasets} 
    finally, any of the tetrahedral datasets could have been modified to obtain polyhedral meshes, but we observed that aggregating elements from $\Dtetuni$, $\Dtetani$ or $\Dtetbcl$ would still generate convex elements, not so different from the original ones.
    We therefore chose to consider only $\Dpolypar$, $\Dpolypoisson$ and $\Dpolyrandom$.
\end{itemize}
Comparing to the datasets defined over the unit cube in \cite{da2017high}, we could say that dataset $\Dvororandom$ is analogous to the ``Random'' discretization, dataset $\Dhexuni$ is equivalent to their ``Structured'' meshes and dataset $\Dvorobcl$ can be considered as a particular case of ``CVT'' discretization.

\medskip
In Table~\ref{table:assumptions} we report some considerations about the datasets and the geometrical assumptions from Section~\ref{subsec:assumptions}.
First, all the considered datasets satisfy assumption \textbf{G1}: the only non convex elements can be found in meshes belonging to polyhedral datasets, but even those elements are still star-shaped, being the union of two tetrahedra.
Datasets originating from parallel and random samplings are not guaranteed to satisfy assumption \textbf{G2}, due to their random nature.
The same holds for datasets from Poisson sampling: even if the points are placed at fixed distance apart, they may still fall arbitrary close to each other around the boundary of $\Omega$.
Datasets originating from the anisotropic sampling instead, are guaranteed to systematically violate assumption \textbf{G2}.
Assumption \textbf{G3} is violated only by $\Dvoropoisson$ and $\Dvororandom$, as it is not possible to bound the number of faces of a random Voronoi cell.
Last, datasets originating from uniform sampling and BCL are guaranteed to satisfy all the geometrical assumptions.
We recall from Section~\ref{subsec:indicator} that to ensure the optimal behavior of the method either assumptions \textbf{G1} and \textbf{G2} or assumptions \textbf{G1} and \textbf{G3} need to be satisfied.
Therefore we are allowed to expect optimal convergence rates over all datasets, except for $\Dvoropoisson$ and $\Dvororandom$.


\begin{table}[htbp]
\caption{Summary of the geometrical assumptions violated by each dataset.}
\label{table:assumptions}
\centering
\begin{multicols}{2}
\begin{tabular}{lccc}
\hline\noalign{\smallskip}
dataset & \textbf{G1} & \textbf{G2} & \textbf{G3}\\
\noalign{\smallskip}\hline\noalign{\smallskip}
$\Dtetuni$ & & & \\
$\Dtetani$ & & $\times$ & \\
$\Dtetpar$ & & $\times$ & \\
$\Dtetbcl$ & & &\\
$\Dtetpoisson$ & & $\times$ &\\
$\Dtetrandom$ & & $\times$ &\\
\noalign{\smallskip}\hline
\end{tabular} 

\columnbreak
\begin{tabular}{lccc}
\hline\noalign{\smallskip}
dataset & \textbf{G1} & \textbf{G2} & \textbf{G3}\\
\noalign{\smallskip}\hline\noalign{\smallskip}
$\Dhexuni$ & & & \\
$\Dhexani$ & & $\times$ & \\
$\Dhexpar$ & & $\times$ & \\
\noalign{\smallskip}\hline
$\Dvorobcl$ & & &\\
$\Dvoropoisson$ & & $\times$ & $\times$ \\
$\Dvororandom$ & & $\times$ & $\times$ \\
\noalign{\smallskip}\hline
$\Dpolypar$ & & $\times$ & \\
$\Dpolypoisson$ & & $\times$ & \\
$\Dpolyrandom$ & & $\times$ & \\
\noalign{\smallskip}\hline
\end{tabular} 
\end{multicols}
\end{table}

\section{Correlations between the quality and the performance}
\label{sec:results}
In order to test the accuracy of the quality indicator $\varrho$ defined in Section~\ref{subsec:indicator}, we evaluate it over each mesh of each dataset from Section~\ref{subsec:dataset:datasets}.
We recall that for an ideal dataset made by meshes containing only equilateral tetrahedra, $\varrho$ would be constantly equal to 1.
We assume this value as a reference for the other datasets: 
the closer is $\varrho$ on a dataset to the line $y~=~1$, the smaller is the approximation error that we expect that dataset to produce.
Moreover, the $\varrho$ slope is indicative of the convergence rate.
Since an ideal dataset would produce an horizontal line, the more negative is the slope, the worse is the convergence rate that we expect.
A positive trend instead, should indicate a convergence rate higher than the one obtained with an equilateral tetrahedral mesh, that is, higher than the theoretical estimates.
This phenomenon is commonly called \textit{superconvergence}.

\medskip
Then, we solve the discrete Poisson problem \eqref{eq:varform} with the VEM \eqref{eq:varform:VEM} described in Section~\ref{sec:vem}, using as groundtruth the function 
\begin{equation}
    \label{eq:groundtruth}
    u(x, y, z) = x^3 y^2 z + x \sin(2\pi xy) \sin(2\pi yz) \sin(2\pi z), \hspace{0.1cm} (x,y,z)\in\Omega.
\end{equation}
For each class of datasets (tetrahedral, hexahedral, Voronoi and polyhedral) we plot the relative $\HONE$-seminorm and $\LTWO$-norm as defined in \eqref{eq:error:H1}, \eqref{eq:error:L2}:
\begin{align*}
    \norm{u-u_h}{0,\Omega}/\norm{u}{0,\Omega}, \qquad
    \snorm{u-u_h}{1,\Omega}/\snorm{u}{1,\Omega},
\end{align*}
and the relative $\LINF$-norm 
\begin{align*}
    \norm{u-u_h}{\infty}/\norm{u}{\infty}, \ \mbox{where }
    \norm{u}{\infty}=\mbox{ess sup}_{x\in\Omega}\snorm{u(x)}{},
\end{align*}
of the approximation error $u-u_h$ as the number of vertices increases (which corresponds to the number of degrees of freedom in our formulation).
The optimal convergence rate of the method, provided by the estimates \eqref{eq:source:problem:H1:error:bound} and \eqref{eq:source:problem:L2:error:bound}, is indicated by the slope of the reference triangle.
In the case of the $\LINF$-norm we do not have such theoretical results.

\medskip
The exercise we propose is to first analyse the values of $\varrho$ on a dataset, computed before solving the problem, and make some predictions on the behaviour of the VEM over it in terms of convergence rate and error magnitude.
Then, looking at the approximation errors actually produced by that dataset, search for correspondences between $\varrho$ and the errors, checking the accuracy of the prediction.
%
%
Clearly, as $\varrho$ does not depend on the type of norm used, we will compare it to an average of the plots for the different norms.

\subsection{Tetrahedral datasets}
\begin{figure}[htbp]
\centering
    \includegraphics[width=.24\linewidth]{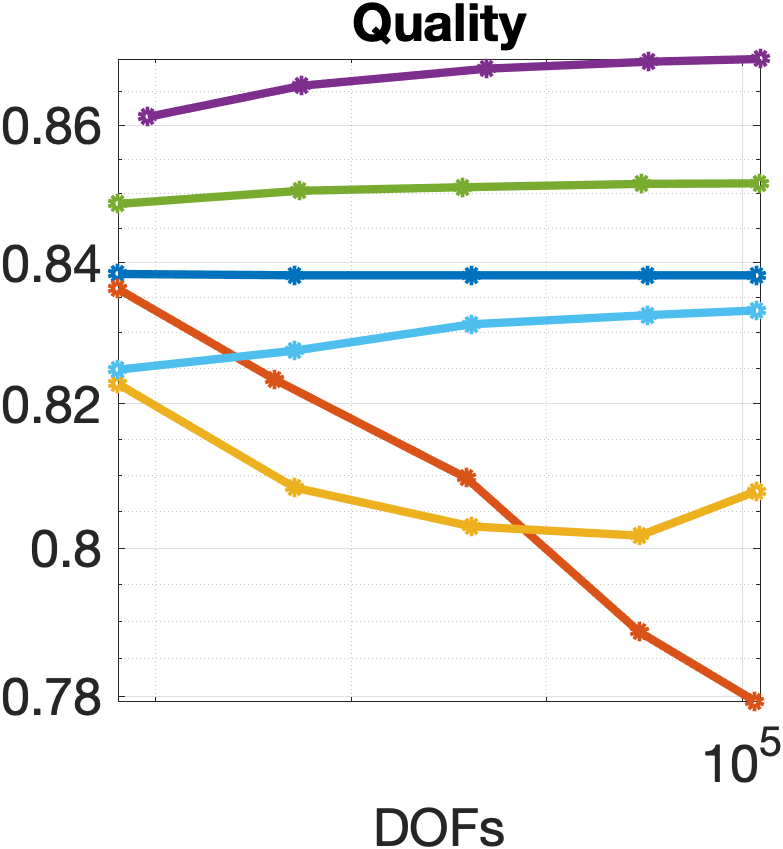}
    \includegraphics[width=.24\linewidth]{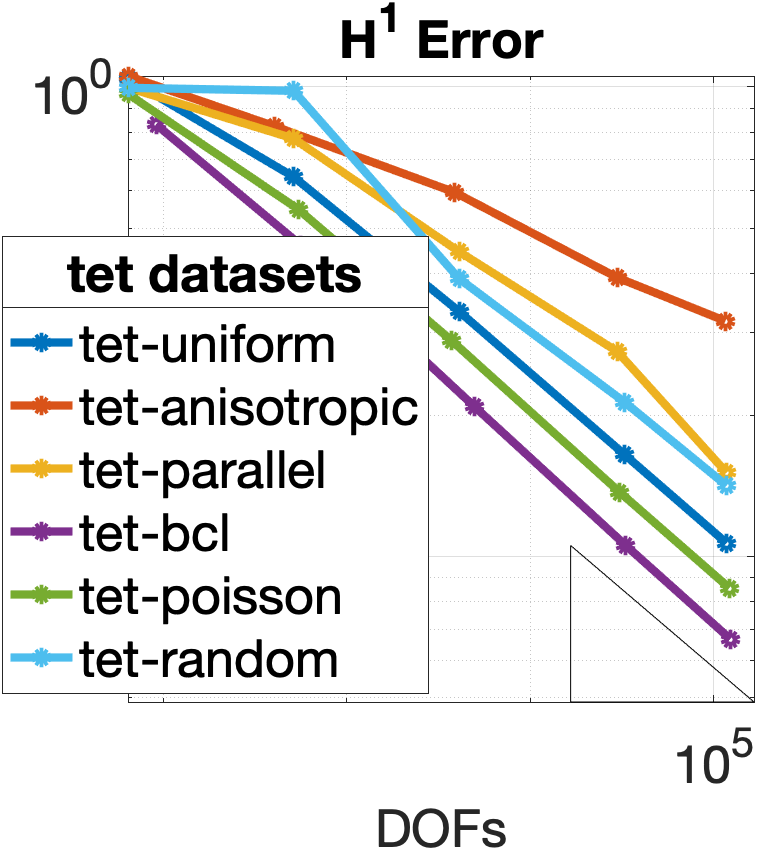}
    \includegraphics[width=.24\linewidth]{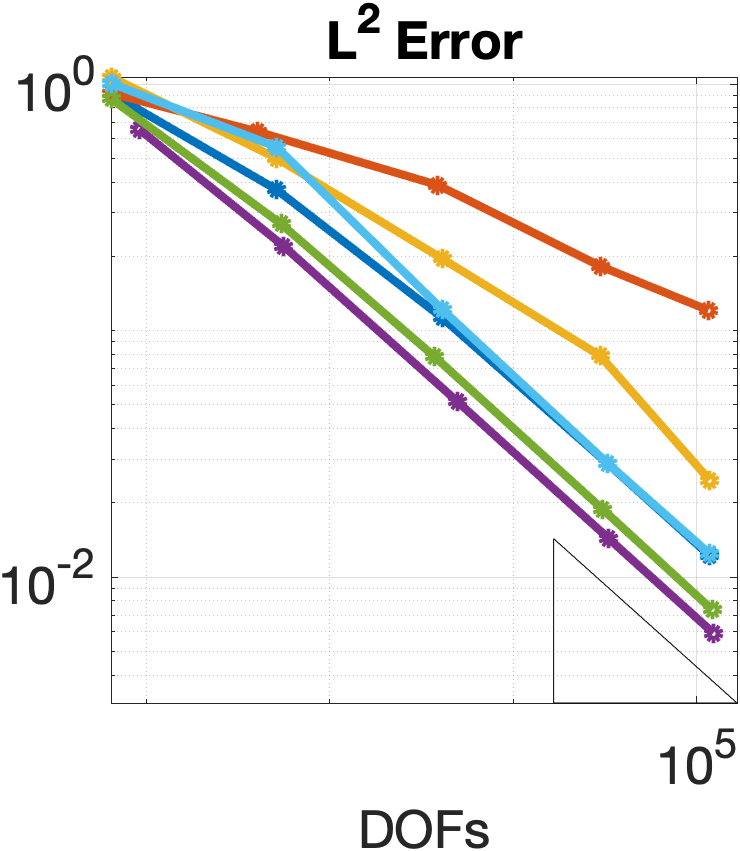}
    \includegraphics[width=.24\linewidth]{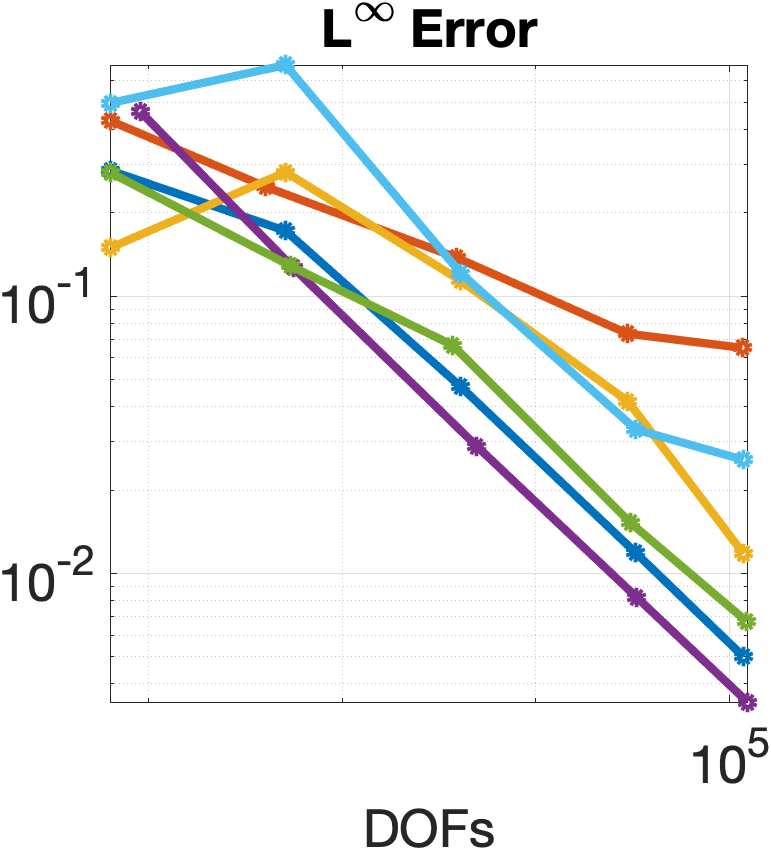}
    \caption{Mesh quality indicator $\varrho$, $H^1$-seminorm, $L^2$-norm and $\LINF$-norm of the approximation errors relative to the tetrahedral datasets.}
\label{fig:tet}
\end{figure}

\paragraph*{Predictions}
Looking at the quality plot in Figure~\ref{fig:tet} (the leftmost) we would say that the VEM should converge with the optimal rate over almost all datasets, as their $\varrho$ tend to get horizontal.
One exception is $\Dtetani$, which has a negative trend and therefore is not expected to converge properly.
We can also observe how the slope of $\Dtetpar$ becomes positive in the last mesh: this should indicate a more than optimal convergence rate.
Regarding the error magnitude, represented by the overall distance from the top of the plot, we can predict that $\Dtetbcl$ will produce the smallest errors, being the one with the highest quality.
This is reasonable because this dataset is composed mainly of equilateral tetrahedra.
We can then order decreasingly the other datasets according to their quality: $\Dtetpoisson$, $\Dtetuni$, $\Dtetrandom$, $\Dtetpar$ and $\Dtetani$, and we expect the errors magnitudes to behave accordingly.

\paragraph*{Verification}
The convergence rates of the tetrahedral datasets in the error plots of Figure~\ref{fig:tet} faithfully respect all the above considerations.
In both $\HONE$ and $\LTWO$-norms, the method converges with the optimal rate (the one suggested by the reference triangle) over all datasets, except for $\Dtetani$, and datasets $\Dtetpar$ superconverges in the last mesh.
We checked the condition numbers for the linear systems that are built on the meshes of $\Dtetani$ and we verified that their values are reasonably small, i.e., in the range $[1,10^6]$.
These values are comparable to the condition numbers seen in the other datasets.
In the $\LINF$-norm the situation is similar, even if $\Dtetrandom$ has an unexpected peak in the last mesh.
The error magnitudes, i.e. the distance of the line from the $y-$axis, perfectly follow the ordering suggested by the quality plot.
The dataset which produces the smallest errors is $\Dtetbcl$.
After that, in the $\HONE$ and $\LTWO$ plots we have $\Dtetpoisson$, $\Dtetuni$ and $\Dtetrandom$, which tend to become very close, then $\Dtetpar$ and last $\Dtetani$.
The situation slightly changes if we look at the $\LINF$ error: in this case $\Dtetuni$ performs better than $\Dtetpoisson$, probably due to a bunch of poor quality elements which do not particularly affect the overall accuracy of the method.

\subsection{Hexahedral datasets}
\begin{figure}[htbp]
\centering
    \includegraphics[width=.24\linewidth]{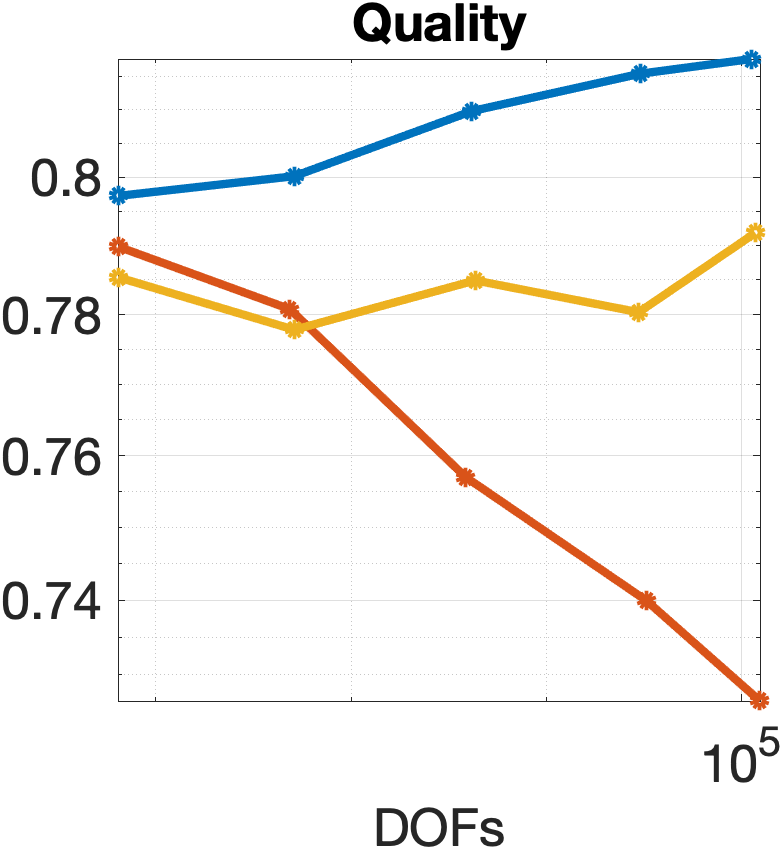}
    \includegraphics[width=.24\linewidth]{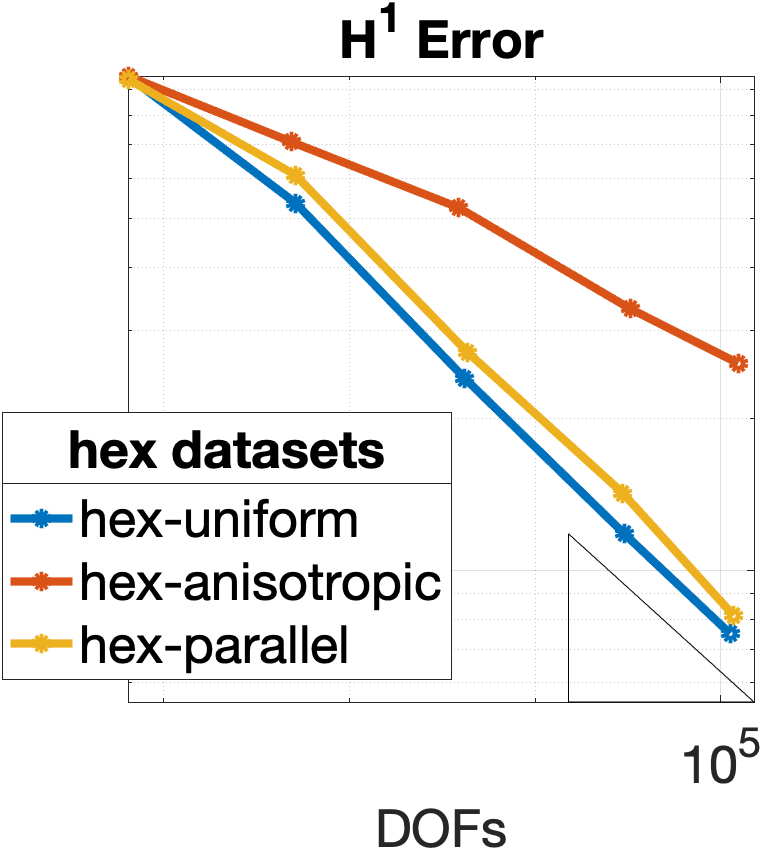}
    \includegraphics[width=.24\linewidth]{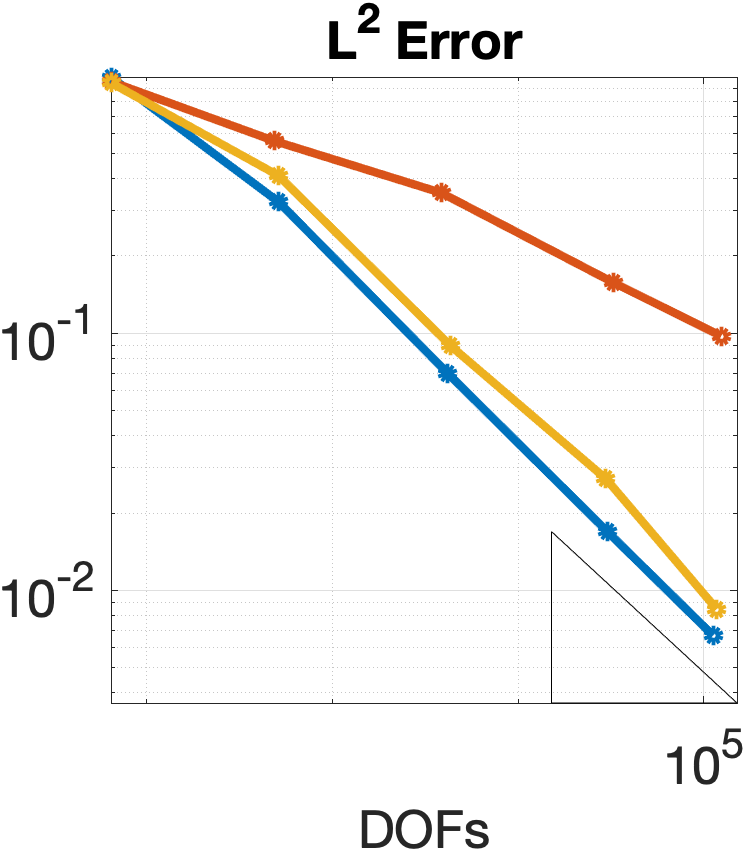}
    \includegraphics[width=.24\linewidth]{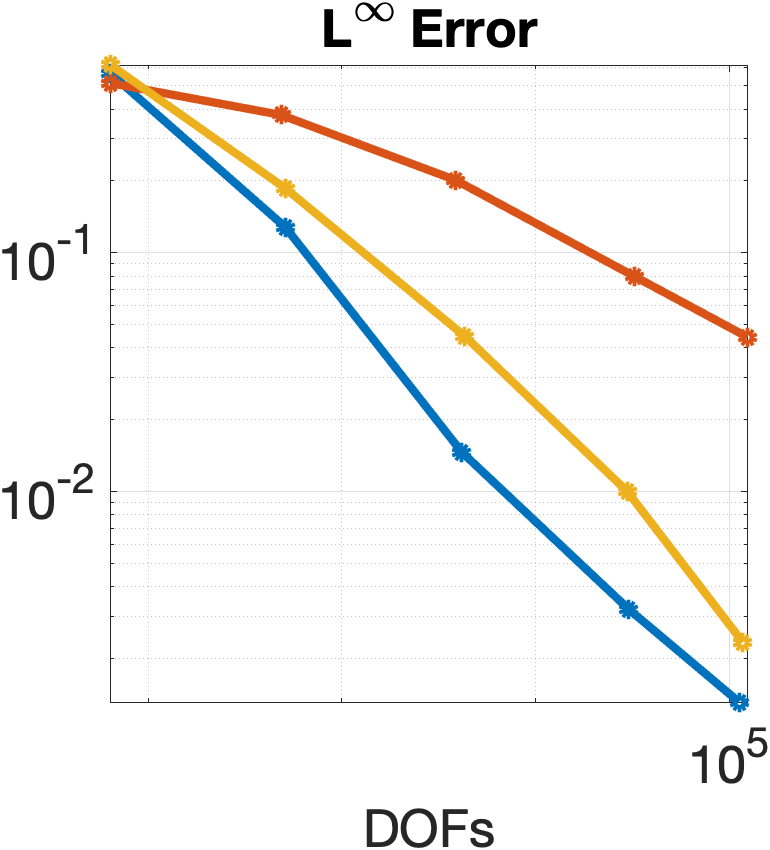}
    \caption{Mesh quality indicator $\varrho$, $H^1$-seminorm, $L^2$-norm and $\LINF$-norm of the approximation errors relative to the hexahedral datasets.}
\label{fig:hex}
\end{figure}

\paragraph*{Predictions}
Results for the hexahedral datasets are shown in Figure~\ref{fig:hex}.
Similarly to what happened for the tetrahedral datasets, the meshes produced by the anisotropic sampling have very poor quality.
While $\Dhexuni$ and $\Dhexpar$ tend to flatten, with the second one increasing in the last refinement, the $\varrho$ value for the meshes of $\Dhexani$ keeps decreasing.
Our prediction is therefore to have optimal convergence on $\Dhexuni$ and $\Dhexpar$ and bad results with $\Dhexani$.
In addition, $\Dhexuni$ is expected to produce smaller errors than $\Dhexpar$.

\paragraph*{Verification}
In the error plots of Figure~\ref{fig:hex} all the predictions are confirmed.
The VEM converges perfectly over $\Dhexuni$ and $\Dhexpar$, with the second one producing higher errors than the first one and improving its convergence rate in the last refinement.
Instead, $\Dhexani$ does not produce a correct convergence rate in the $\HONE$ and $\LTWO$ plots, and also in the $\LINF$ plot exhibits a significantly slower rate with respect to the other datasets.
Also in this case, the condition numbers for the linear systems are not particularly bigger than the ones of the other datasets.

\subsection{Voronoi datasets}
\begin{figure}[htbp]
\centering
    \includegraphics[width=.24\linewidth]{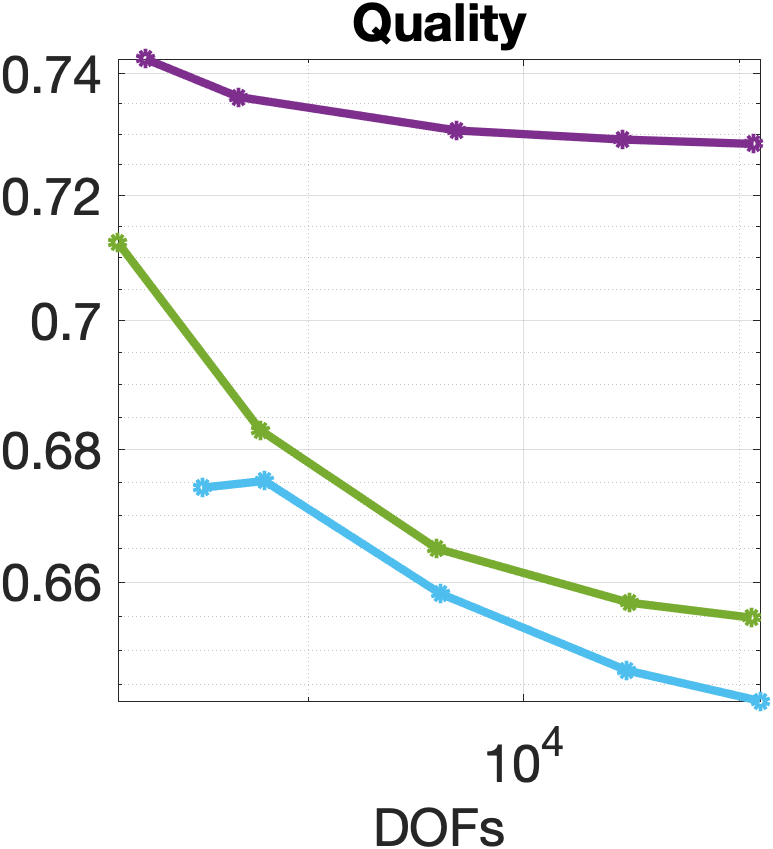}
    \includegraphics[width=.24\linewidth]{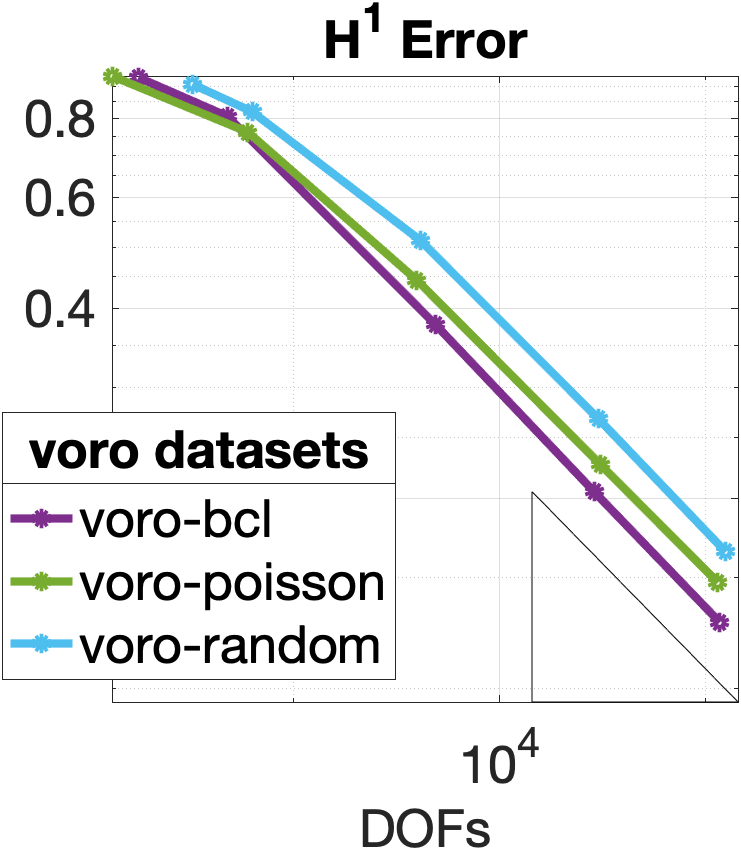}
    \includegraphics[width=.24\linewidth]{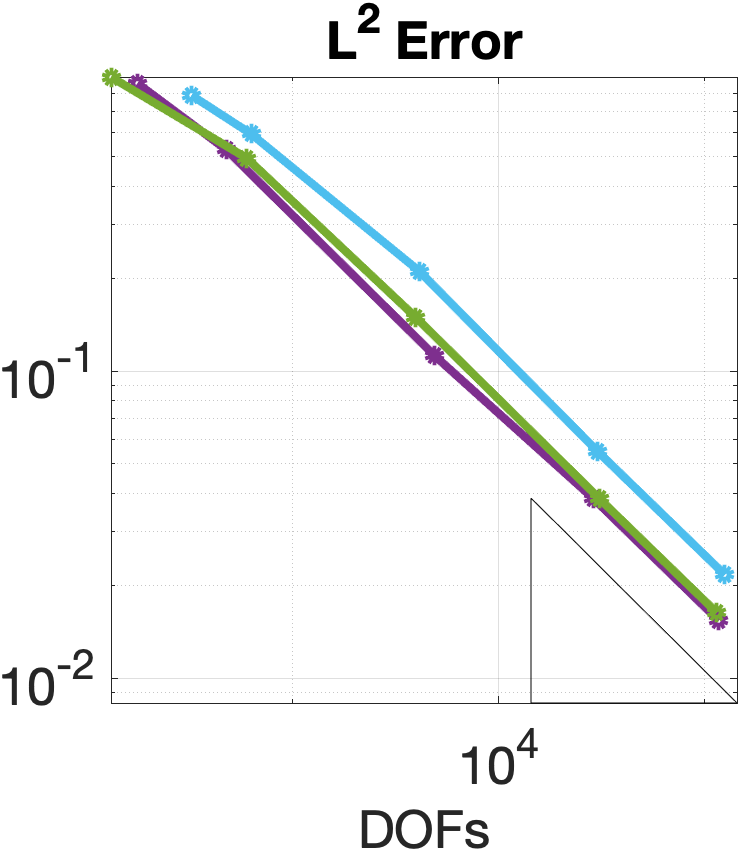}
    \includegraphics[width=.24\linewidth]{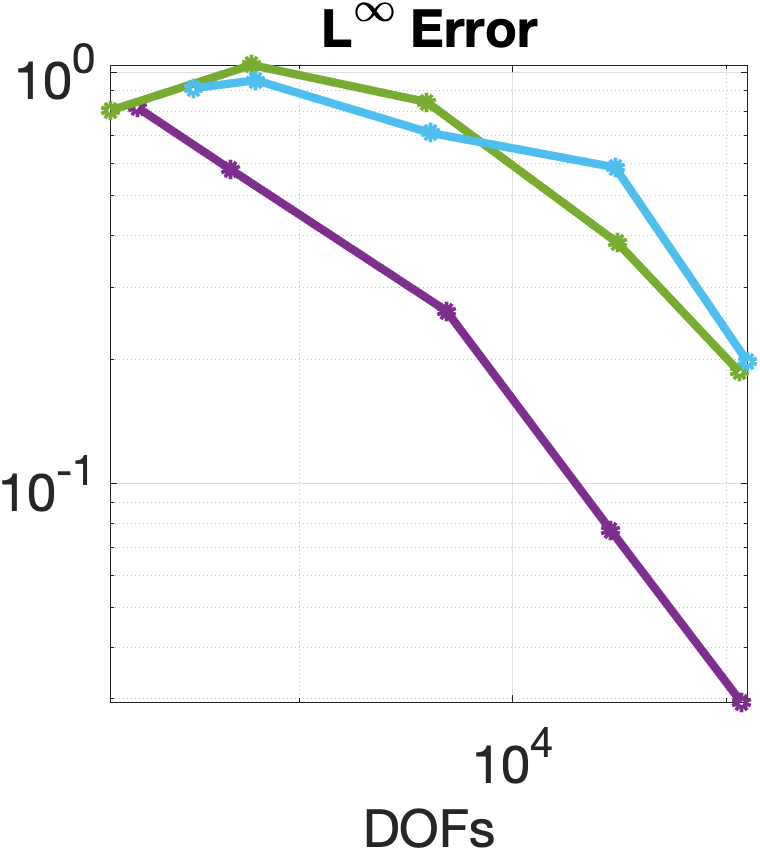}
    \caption{Mesh quality indicator $\varrho$, $H^1$-seminorm, $L^2$-norm and $\LINF$-norm of the approximation errors relative to the Voronoi datasets.}
\label{fig:voro}
\end{figure}

\paragraph*{Predictions}
In Figure~\ref{fig:voro}, results relative to the Voronoi datasets are shown.
The quality of all three datasets tend to stabilize to a constant value, and this makes us presume a correct convergence rate for all of them.
We can expect $\Dvorobcl$ to produce smaller errors than the other two, and $\Dvororandom$ to be the less accurate.

\paragraph*{Verification}
Looking at the $\HONE$ and $\LTWO$ error plots we notice how all datasets converge properly, and the accuracy of the approximation follows the order foreseen by the indicator: $\Dvorobcl$, $\Dvoropoisson$ and $\Dvororandom$.
The $\LINF$ plot is less similar to the other two in this case, but still we can recognise a common pattern.

\subsection{Polyhedral datasets}
\begin{figure}[htbp]
\centering
    \includegraphics[width=.24\linewidth]{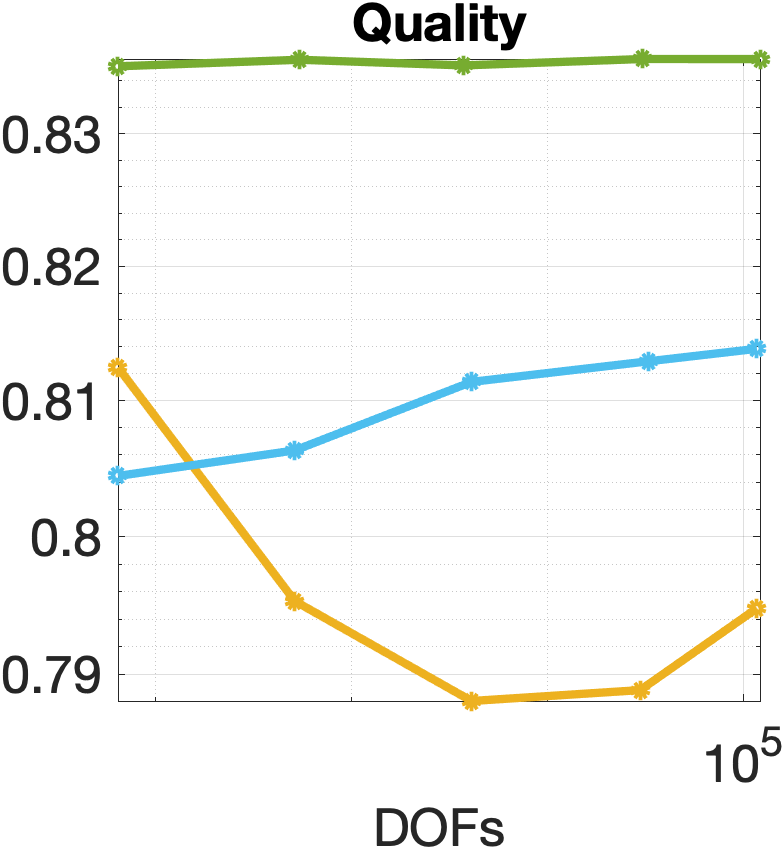}
    \includegraphics[width=.24\linewidth]{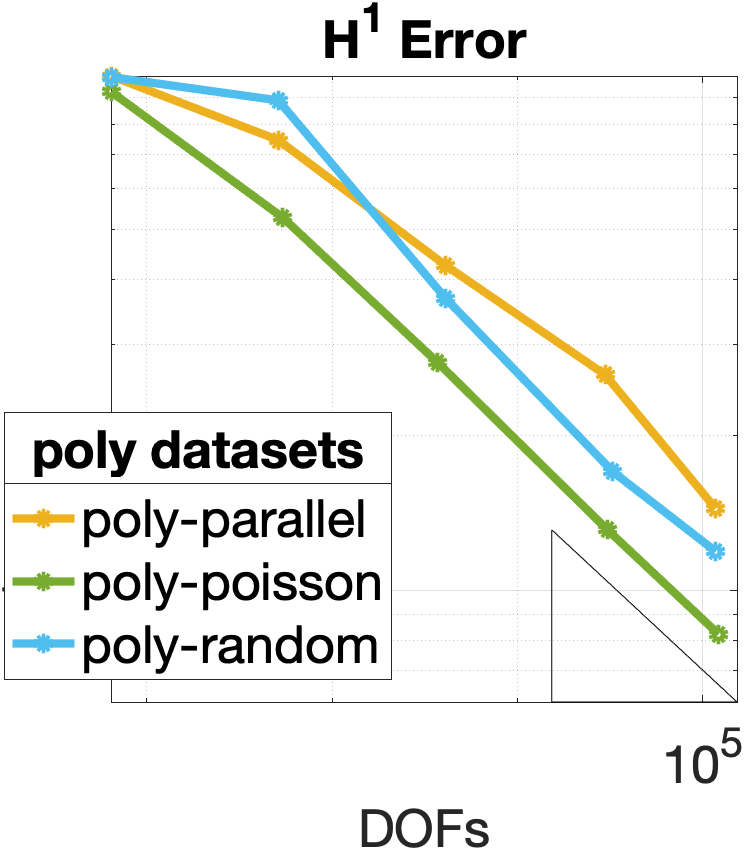}
    \includegraphics[width=.24\linewidth]{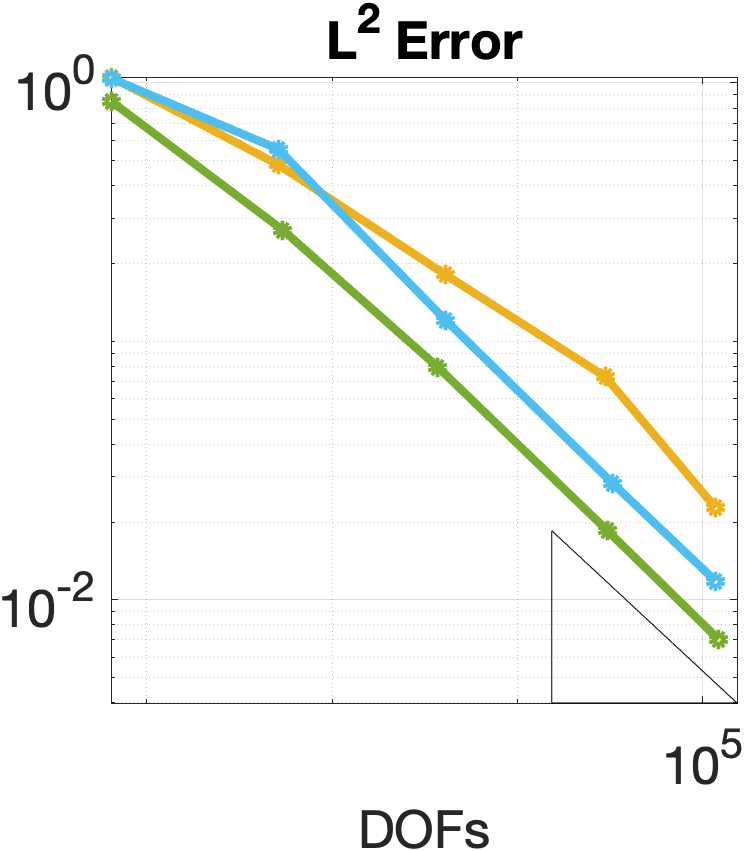}
    \includegraphics[width=.24\linewidth]{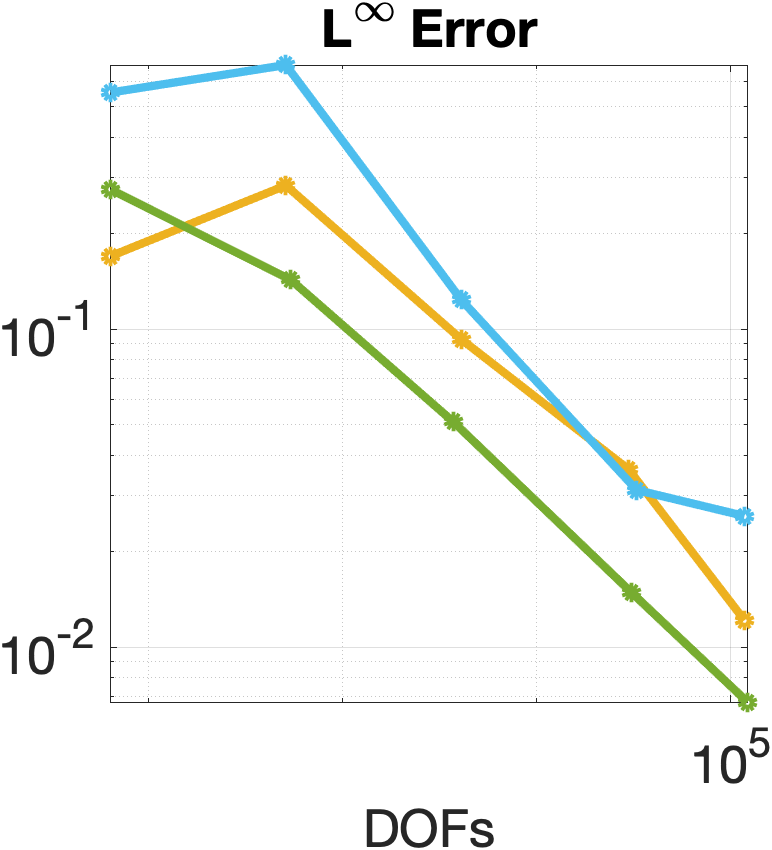}
    \caption{Mesh quality indicator $\varrho$, $H^1$-seminorm, $L^2$-norm and $\LINF$-norm of the approximation errors relative to the polyhedral datasets.}
\label{fig:poly}
\end{figure}

\paragraph*{Predictions}
Last, in Figure~\ref{fig:poly} we report the analysis of the polyhedral datasets.
The indicator $\varrho$ suggests that $\Dpolypoisson$ and $\Dpolyrandom$ converge perfectly.
Regarding $\Dpolypar$, the indicator seems to flatten and then increases in the last refinement.
The convergence should therefore be optimal for the first meshes and more than optimal for the last one.
The most accurate dataset should be $\Dpolypoisson$ and the least accurate $\Dpolypar$.

\paragraph*{Verification}
The method performs essentially as expected.
All datasets produce optimal rates and $\Dpolypar$ converges even faster than the reference in the last refinement.
Dataset $\Dpolyrandom$ has a peak in the last mesh with the $\LINF$ error: this is similar to what happened with $\Dtetrandom$ and it probably due to the same bad-shaped element (remember that tetrahedral and polyhedral meshes differ only for the $20\%$ of their elements).
Concerning the errors magnitude, as foreseen by the indicator, $\Dpolypoisson$ is the most accurate and then we have $\Dpolyrandom$ and $\Dpolypar$.

\section{Conclusions}
\label{sec:conclusions}
We conclude with some more general considerations on the results
obtained in the previous section.
When all three geometrical assumptions are respected, that is, with
meshes from uniform sampling and BCL, the performance are obviously
optimal.
It is important to notice how $\Dtetani$ and $\Dhexani$ are the only
two cases in which the VEM underperforms, despite
satisfying both assumptions \textbf{G1} and \textbf{G3}.
In these cases, the strong violation of assumption \textbf{G2} significantly impacts on the proper convergence of the method.
In other similar situations, e.g. with parallel samplings, \textbf{G2}
is violated less heavily and the method works as expected.
Vice-versa, $\Dvoropoisson$ and $\Dvororandom$ satisfy \textbf{G1} but
not \textbf{G2} nor \textbf{G3}, and the VEM still manages to converge
on them.
These results are not unexpected, as the geometrical assumptions are only \textit{sufficient} conditions for the convergence of the method, but confirm our suspects that the current restrictions on the meshes are probably more severe than necessary.

The mesh quality indicator has been able to properly predict the behaviour of the VEM, both in terms of convergence rate and error magnitude, up to a certain precision.
It showed up to be particularly accurate when compared to the $\HONE$ and the $\LTWO$-norms of the error, while it not always managed to capture the oscillations of the $\LINF$-norm.

In conclusion, the relationship between the geometrical assumptions
respected by a dataset and the performance of the VEM on it is not so
obvious, especially when we try to violate at least one of them.
In those situations, our mesh quality indicator turns out to be
particularly useful in predicting the result of the numerical
approximation.
Its effectiveness lies in the ability of capturing a qualitative
measure of the violation of the single assumptions.

\medskip
We are currently working on a software library capable of splitting or merging the elements of a mesh in order to maximize an energy functional based on the quality indicator.
This software is capable of spotting the most pathological elements in a mesh (the ones with the poorest quality) and either aggregate them with a neighbor or split them into smaller parts, improving the global quality of the mesh.
We believe this could be extremely useful in VEM simulations, as a higher quality mesh leads to cheaper and more accurate approximations.

As a future work, we plan to further investigate the geometrical
assumptions involved in the three-dimensional VEM analysis, for instance an extension of assumption \GAs{4}, and to study the convergence of the VEM with order
$k>1$.
In \cite{sorgente2021role} it was shown that the behaviour of the VEM
over a polygonal mesh does not drastically change for different values
of $k$, and therefore the quality measured by $\varrho$ can be
considered as a reliable indicator for all the orders of the method.
We suspect this could be true also for the three-dimensional
formulation, but more detailed studies are required on this topic.
Last, it would be interesting to see if the quality indicator could be adapted to other numerical schemes different from the VEM, such as the Discontinuous Galerkin method \cite{cockburn2012discontinuous}, by opportunely defining a new set of scalar functions based on appropriate geometrical assumptions.

\section*{Acknowledgments}
This paper has been realised in the framework of ERC Project CHANGE,
which has received funding from the European Research Council (ERC)
under the European Union's Horizon 2020 research and innovation
program (grant agreement no.~694515).


\bibliographystyle{plain}

\end{document}